\documentclass[12pt,reqno,psfig,openbib]{amsart}
\usepackage{amssymb,amsmath,graphicx,amsfonts,euscript}%,showkeys}
\usepackage{fancyhdr}
\usepackage{epsfig}
\usepackage{setspace}
\usepackage{dsfont}
\usepackage{subcaption}
\usepackage{xcolor}
\usepackage{color}

\newcommand{\e}{\varepsilon}

\newcommand{\f}{\frac}

\newcommand{\be}{\begin{equation}}
\newcommand{\ee}{\end{equation}}

\numberwithin{equation}{section}
\numberwithin{figure}{section}
\numberwithin{table}{section}

\newtheorem{rem}{Remark}[section]

\title[Deep neural network motivated by Legendre-Galerkin approximation]{Deep neural network for solving differential equations motivated by Legendre-Galerkin approximation}

\author[B. Chudomelka]{Bryce Chudomelka}
\address[Bryce Chudomelka]{Department of Mathematics and Statistics, San Diego State University, San Diego, CA, USA.}
\email{bchudomelka@sdsu.edu}

\author[Y. Hong]{Youngjoon Hong}
\address[Youngjoon Hong]{Department of Mathematics and Statistics, San Diego State University, San Diego, CA, USA.}
\email{yhong2@sdsu.edu}

\author[H. Kim]{Hyunwoo Kim}
\address[Hyunwoo Kim]{Department of Computer Science, Korea University, Seoul, Korea.}
\email{hyunwoojkim@korea.ac.kr}

\author[J. Park]{Jinyoung Park}
\address[Jinyoung Park]{Department of Computer Science, Korea University, Seoul, Korea.}
\email{lpmn678@korea.ac.kr}

\date{\today}

\begin{document}
    \maketitle
    
    \begin{abstract}
    Nonlinear differential equations are challenging to solve numerically and are important to understanding the dynamics of many physical systems.
    Deep neural networks have been applied to help alleviate the computational cost that is associated with solving these systems. We explore the performance and accuracy of various neural architectures on both linear and nonlinear differential equations by creating accurate training sets with the spectral element method. 
    Next, we implement a novel Legendre-Galerkin Deep Neural Network (LGNet) algorithm to predict solutions to differential equations.
    By constructing a set of a linear combination of the Legendre basis, we predict the corresponding coefficients, $\alpha_i$ which successfully approximate the solution as a sum of smooth basis functions $u \simeq \sum_{i=0}^{N} \alpha_i \varphi_i$.
    As a computational example, linear and nonlinear models with Dirichlet or Neumann boundary conditions are considered.
    % We propose a residual of the weak formulation of the differential equations for loss functions of the neural network.
\end{abstract}
    \noindent
    \section{Introduction}
    Partial differential equations (PDEs) are often difficult to solve, frequently requiring varying numerical approximations, where each method has its own computational cost and accuracy. The robust history and research into discretization and numerical methods for solutions of PDEs raises an interesting problem for deep neural networks (DNNs), specifically because of the rich mathematical theory and numerical analysis. In this work, we develop a novel DNN for solving differential equations based on Legendre-Galerkin approximation by leveraging accurate solutions in a supervised learning setting to find solutions given a forcing function input, $f$, and boundary conditions. For accurate and reliable training, we utilize a residual of the weak formulation of the differential equations (DE) as a loss function. Different types of equations including fluid and wave models are implemented to demonstrate numerical performance of the DNNs.
    
    Deep learning is the process by which NNs learn underlying patterns to robust data sets through affine transformations and can be generalized through a process known as cross-validation; these networks are \textit{deep} because they have many hidden layers \cite{Allen74, arora2016understanding, Goodfellow-et-al-2016}. NNs have demonstrated excellent utility in the approximation of continuous functions and have achieved state-of-the-art performance in many fields, while \textit{deep} neural networks have shown performance as a function of depth, \textit{e.g.}, adding more layers to the network improves accuracy \cite{Bar_Sinai_2019, Cybenko1989ApproximationBS, he2015deep, HORNIK1989359, Goodfellow-et-al-2016, Stinchcombe1989, zhuang2020learned}. Recently, mathematical and numerical approaches of DNNs have been studied in computer vision to improve adversarial robustness \cite{kim2020sspnet}. This research is loosely related, developing a deeper understanding of the implementation of DNNs to numerical solutions.  The central goal of our research is to investigate an accurate numerical solution which can be approximated and predicted by a DNN.
    
    Previous works have demonstrated the success of neural architectures when applied to solutions of differential equations \cite{Bar_Sinai_2019, kharazmi2019variational, Raissi1026, RUDD2015277, zhuang2020learned}. 
    The authors in \cite{Bar_Sinai_2019, zhuang2020learned} introduce data-driven discretization, a method for learning optimized approximations to PDEs based on traditional finite-volume (or difference) schemes. The algorithm uses neural networks to estimate spatial derivatives, which are optimized end to end to best satisfy the equations on a low-resolution grid.
    Recently, the Physics Informed Neural Networks (PINN) introduces a novel methodology for finding solutions to complex dynamical systems utilizing automatic differentiation; see e.g.  \cite{Raissi1026, RAISSI_PINN, kharazmi2019variational, RK18} among many other references. 
    In particular, \cite{kharazmi2019variational} utilizes the variational form, \textit{i.e.}, the weak form of the derivative, to increase the accuracy of their network.

    Our research differs from the previous work by exploring different neural architectures to obtain an accurate solution that we can find. 
    We use a deep convolutional neural network (CNN) to achieve \textit{sufficiently} accurate solutions by predicting coefficients of spectral approximation based on Legendre-Galerkin methods \cite{GO77, STW11}.  
    This approach is partially related to the data-driven discretization in \cite{Bar_Sinai_2019, zhuang2020learned}, but their neural networks predict coefficients of numerical derivatives such as finite difference methods.
    In addition, accuracy of the training data set generated by the finite volume methods is limited even if high resolution methods are implemented. 
    When fitting a model, we introduce a residual of numerical integration inherited from a weak formulation of the differential equations.
    This choice is natural as the architecture mainly relies on spectral element methods constructed by the weak formulation of the differential equations. 
    The choice of loss function is close to \cite{kharazmi2019variational}, but was done independently without prior knowledge and utilizes different neural architectures, as well as a novel algorithm; the Legendre-Galerkin spectral element method (LGSEM) with DNNs. 
    Hence, the test function is a linear combination of Legendre polynomial which satisfy boundary conditions of the model.
    It is noteworthy that the methodology under consideration can incorporate various choices of polynomial basis functions into the novel neural network architecture. 
    More precisely, there are many feasible choices of basis functions, such as Fourier series, Chebyshev polynomials, or Jacobi polynomials.
    Hence, our architecture is flexible and extendable to other numerical approximation. 
    In addition, since the coefficients of the spectral approximation is predicted by the network, the main structure of the numerical approximation with the polynomial basis is maintained.
    In this regard, various existing numerical methods such as enriched space method are applicable to the proposed network architecture by adding a proper basis function; for more details, see e.g. \cite{HJT14, CHT20}.
   
    The methodology of constructing the test data set, as well as the architecture of the network, used to predict the solution should play a critical component in the accuracy of the predicted solutions. A network that is trained on noisy or corrupt data could suffer in accuracy, especially in the application of neural networks to nonlinear DEs. 
    A training data set generated by that is pristine offers a training data set that is highly accurate, via the Spectral Element Method (SEM), for a neural network to learn the governing dynamics of a DE. The SEM can achieve spectral accuracy with a relatively small number of collocation points, which places the accuracy of our method on the order of machine precision. Thus, our data set will be as accurate as possible in order for the network to predict the dynamics.
    
    In the calculations presented here, the model consists of convolutional layers connected with a nonlinear activation function between each layer. CNNs are used in a wide range of tasks such as facial recognition, object detection, semantic segmentation, natural language processing, and recently, time-series data, for state-of-the-art results \cite{cui2016multiscale, kim2014convolutional, AlexNet, ogunmolu2016nonlinear, tang2020rethinking}. CNNs are a series of linear operations, and in the context of a DNN, can be layered together with nonlinear operations, proving useful for a universal function approximator \cite{Cybenko1989ApproximationBS, HORNIK1989359}. We train multiple configurations of CNN neural architectures with different nonlinear operations, known as activation functions, on robust data sets generated with the SEM to perform supervised learning. Each configuration is then tested using an out-of-sample set of 1000 randomly sampled solutions that are not present in the input data set. This allows us to measure the ability of our network to generalize to unknown data.

    In this paper, we develop a version of CNNs based on the Legendre-Galerkin framework to find a numerical solution of DEs via a deep neural network. The networks predicts coefficients, $\alpha_i$ of spectral approximation consisting of Legendre polynomial basis $\varphi_i$ such that $u \simeq \sum_{i=0}^N \alpha_i \varphi_i$, where $u$ is a solution of the DE. In addition, for accurate and reliable training, we utilize a residual of weak formulation of the DEs as a loss function. Different types of equations including fluid and wave models are implemented.
    
    The article is organized as follows: In Section 2 we address the nonlinear function approximation of DNNs, and introduce a novel architecture of the DNNs based on numerical approximation. In Sections 3 and 4, linear and nonlinear models equipped with different boundary conditions are implemented. A sequence of numerical experiments are presented to demonstrate the numerical performance of the proposed method. We conclude the paper with a summary section 5.

    \section{Data-Driven Numerical Methods}
    Consider a set of differential equations:
    \be
    \begin{split}
    	& \mathcal{F}(u, u_x, u_{xx}) = f(x), \quad x \in [-1,1],\\
    	& \mathcal{B}(u, u_x) = g(x), \quad\quad \text{at} \quad  x=-1,1,
    	\label{eq:PDE}
    \end{split}
    \ee
    where $\mathcal{F}$ is a linear or nonlinear operator and $\mathcal{B}$ is a boundary operator. We begin by posing the question, ``Can a neural network learn a {sufficiently} accurate solution to a differential equation when given the forcing function?"
    
    \begin{figure}
        \begin{center}
            \includegraphics[width=0.9\textwidth]{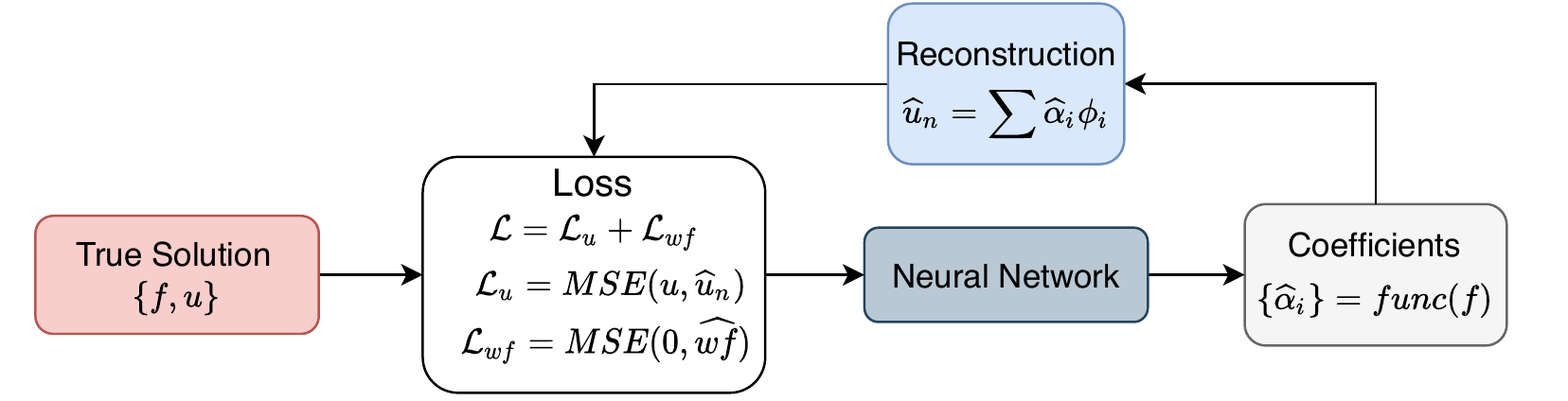}
            \caption{The Legendre-Galerkin Deep Neural Network Algorithm. This diagram demonstrates the training process for the LGNet. The network's input, the forcing function $f$, and the output, a set of coefficients $\{\widehat{\alpha}_i\}$, that are used to reconstruct the predicted solution, $\widehat{u}$. The next step in the training loop is to measure the loss, $\mathcal{L}$, between the ground truth solution $u$ and the predicted solution. The weak form loss is computed by minimizing the weak form for the differential equation. The iterative process then repeats until finished.}
            \label{fig:architecture}
        \end{center}
    \end{figure}
    
    The heart of the algorithm behind our method is to make use of a version of spectral methods, LGSEM, which is one of the compelling methods in numerical analysis, to create robust training sets. One of the major advantages of the LGSEM is to achieve spectral accuracy in a numerical solution \cite{STW11, GO77}. In other words, with a relatively small number of numerical basis, a Legendre basis, the numerical errors are on the order of machine precision ($10^{-14}\sim 10^{-16}$). Thus, guaranteeing that our input data set is pristine so that we can obtain accurate solutions to Equation \eqref{eq:PDE} with a NN model.
    
    The spectral methods rely on a global discretization. The traditional way to introduce them starts by approximating the function as a sum of very smooth basis functions:
    \be \label{eq:SEM}
    	u(x) \simeq \sum_{k=0}^{N-1} \alpha_k \phi_k(x),
    \ee
    where $N$ is the number of the global basis function. In practice, there are many feasible choices of basis functions, such as Fourier series, $\exp(i k x)$, Chebyshev polynomial, $T_k(x)$, or Legendre polynomial, $L_k(x)$.

    In this paper, we focus on the Legendre polynomials.
    The main advantage of Legendre polynomials is that they are mutually orthogonal in the standard $L^2$-inner product, so the analysis of Legendre spectral methods is relatively simple. 
    In addition, an efficient algorithm, which lead to systems with sparse matrices for the discrete variational formulations, is available by constructing appropriate base functions.
    These properties are useful when we derive loss functions based on the weak formulation later. 
    The LGSEM uses a linear combination of Legendre polynomial basis functions which allows our method to adapt to different boundary conditions including Dirichlet, Neumann, and Robin boundary conditions.
    
    Let us define the Legendre polynomial, $L_k(x)$, on the Gauss-Lobatto quadrature, which plays an important role on numerical differentiation and integration. The collocation points form a non-uniform discretization in the spatial domain, but leads to an accurate numerical integration and differentiation. 
    Hence, we set the global basis functions
    \be
    	\phi_k(x) := L_k(x) + a_k L_{k+1} + b_k L_{k+2}(x), \quad x \in [-1,1],
    	\label{eq:global_basis_func}
    \ee
    where $a_k$ and $b_k$ are determined by the boundary conditions of the underlying differential equations. 
    Applying the SEM, we can obtain an accurate numerical solution given by \eqref{eq:SEM}; for more details on the LGSEM, see e.g. \cite{STW11}.
    The set of numerical solutions generated by the SEM is used for the training data set in our neural network.
    
    For any given $f$, we propose a Legendre-Galerkin Deep Neural Network (LGNet) algorithm to find the coefficients, $\{\widehat{\alpha}_k\}$, from which we can reconstruct the corresponding predicted numerical solution, $\widehat{u}$, to the known numerical solution $u(x)$. 
    If our random input data is given by
    \be
    	f(x) = m_1 \sin(\pi w_1 x) + m_2 \cos(\pi w_2 x),
    	\label{eq:forcing}
    \ee
    where parameters, $m_1, m_2, w_1, w_2$, follow normal or uniform distributions, then the corresponding output will be
    \be
    	\{\widehat{\alpha}_k\} \implies \widehat{u} = \sum^{N-1}_{k=0} \widehat{\alpha}_k \phi_k,
    \ee
    where $\{\widehat{\alpha}_i\}$ describes the set of learned coefficients that are used to reconstruct the predicted numerical solution, $\widehat{u}_N$.
    
    We can generate high quality data sets, with an arbitrary number of solutions, using the SEM. These solutions are then fed through the LGNet algorithm, reconstructed, and then the difference between the predicted solution and the actual solution are measured as the primary metric of performance; see Figure \ref{fig:architecture} for reference. This can be done using either the mean absolute error defined as,
    \be
        \label{eq:mae}
        MAE\{g,\widehat{g}\}:=\frac{1}{n}\sum_{i=0}^N\left|g_i-\widehat{g}_i\right|,
    \ee
    or using the mean squared error defined as,
    \be
        \label{eq:mse}
        MSE\{g,\widehat{g}\}:=\frac{1}{n}\sum_{i=0}^N\left(g_i-\widehat{g}_i\right)^2,
    \ee
    for $n$ arbitrary functions $g_i$, and its predicted function $\widehat{g}_i$. Other loss metrics can be used, such as the root mean square error or relative $\ell^2$ error, but we found \eqref{eq:mae} and \eqref{eq:mse} to outperform others. It was observed that the mean square error provides a better metric for generalization of a neural network with regards to the accuracy of the predicted solution. Previous work of \cite{Bar_Sinai_2019} utilize \eqref{eq:mae} in their loss metric, but we observed that \eqref{eq:mse} give the best performance for the LGNet. The predicted coefficients for the global basis vectors, $\{\widehat{\alpha}_i\}$, were more accurate for \eqref{eq:mse} than \eqref{eq:mae}.  We found this to be a consequence of the sensitivity of reconstructing accurate solutions to the learned coefficients when using the SEM. Although the coefficients are not used in the loss metric, we did have access to them for reference indicating that the choice of basis is arbitrary and can be generalized to improve scalability.

    Networks can be further improved by introducing proper loss functions, which act as normalization factors and were also previously used in related works \cite{andrew2007scalable, Bar_Sinai_2019, RAISSI_PINN}. In this regard, we introduce a {\it weak formulation} of \eqref{eq:PDE}. We found that measuring the loss metric in Equation \eqref{eq:mae}, or \eqref{eq:mse}, on the weak form improved the accuracy of our networks when coupled with the loss of the solution. Heuristically, we search for the optimal neural architecture which yields the best performance via a predefined metric; the loss function, $\mathcal{L}$, is the objective function we wish to minimize to gauge performance of the network, but we measure the true performance of a network on its \textit{mean} relative $\ell^2$ error over a test set. Once the structure of the network has been determined then experiments must be carried out to determine the optimal parameters of that architecture. We are interested in finding a neural architecture which gives us the most accurate solution, which leads us to measure the loss between predicted solutions, $\widehat{u}$, and the ground truth solutions contained in the training data set, $u$, as well as using the weak form. The performance of the network is further improved by avoiding the high order numerical derivatives. By using integration by parts on high order derivatives, we circumvent high order numerical differentiation that introduce large numerical errors. This adds another contribution to our network performance metric that acts as a regularization component. Similar observation has been made in previous works \cite{RAISSI_PINN, Raissi1026}, although our work was done independently without prior knowledge.
    
\begin{rem}
Other polynomials such as Chebyshev or Jacobi polynomials can be combined with the proposed DNN architecture since Chebyshev spectral galerkin methods or Jacobi spectral galerkin method have been well developed by those polynomial basis functions. 
\end{rem}

    We now construct the deep neural network based on the convolutional neural network. 
    In \cite{Bar_Sinai_2019}, the authors introduced a CNN with finite volume discretization. However, the method is limited to the periodic boundary conditions, and the accuracy should be similar or smaller than their training sets, which is generated by the 2$^{nd}$ or 3$^{rd}$ order numerical methods. However, our algorithm works for various boundary conditions including Dirichlet, Neumann, and Robin boundary conditions, as long as we choose a right set of Legendre polynomial. We note that the choice of $a_k$ and $b_k$ for the Legendre basis function in \eqref{eq:global_basis_func} gives an exact boundary condition. After predicting the coefficients $\alpha_k$, the numerical solution obtained from reconstruction possesses an exact boundary condition. Hence, numerical errors from the boundary conditions are very small. In addition, since our methods is based on the spectral element methods, we expect small numerical errors.

    \section{Linear model}
    \begin{figure}
        \begin{subfigure}[t]{.45\textwidth}
            \begin{center}
                \includegraphics[width=0.4\textwidth]{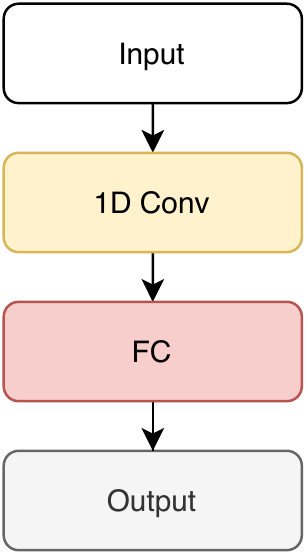}
                \caption{\texttt{Linear} Architecture}
                \label{subfig:architecture_linear}
            \end{center}
        \end{subfigure}
        \hspace{0.5cm}
        \begin{subfigure}[t]{.45\textwidth}
            \begin{center}
                \hspace{1cm}
                \includegraphics[trim=0.0cm 1.5cm 1.85cm 0cm, clip, height=0.25\textheight]{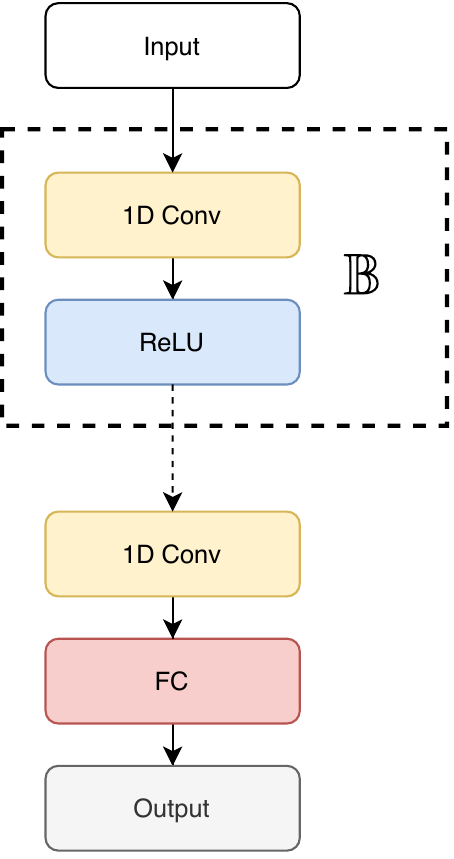}
                \caption{\texttt{NetA} Architecture}
                \label{subfig:architecture_netA}
            \end{center}
        \end{subfigure}
        \caption{The \texttt{NetA} neural network architecture is comprised of multiple blocks, $\mathbb{B}$, that add depth to the network. This diagram depicts a \texttt{NetA} architecture with just 1 block that is comprised of a 1D-convolutional operation with $\mathbb{F}$ filters followed by a \texttt{ReLU} operation. The final convolution in the network does not have the \texttt{ReLU} operation and is then flattened into a fully connected network. The output of this network is a set of coefficients, $\{\widehat{\alpha}_i\}$ that are then used to reconstruct the solution.}
        \label{fig:linear}
    \end{figure}
    
    We want to observe the effectiveness of how a NN can be generalized to find solutions to differential equations. It is most appropriate to begin with an example that is not terribly difficult, yet still non-trivial, to validate the efficacy of this approach. 
    The convection diffusion equation (CDE) is a challenging problem, we refer to this as our paradigm problem, for us to investigate since the structure of the equation can be generalized to many other interesting problems such as Burgers' equations; see e.g. Section \ref{s:4}. All of our experiments were carried out on an Intel Core i9-10900K with an NVIDIA GeForce RTX 2080 SUPER GPU. \\
        
    \subsection{Convection Diffusion Equation}~\\
    
        To deliver our idea, we start with a paradigm problem:
        \be
        \label{eq:paradigm}
        \begin{cases}
        	 -\e u_{xx} - u_x = f,\\
        	 u(-1) = u(1) = 0,
        \end{cases}
        \ee
        where $\e = O(1)$ is a diffusion parameter. The paradigm problem poses an interesting challenge for a neural network. Once the paradigm problem is solved, our algorithm can naturally be generalized to nonlinear, time-dependent, or 2D problems. Notice that homogeneous Dirichlet boundary conditions are used but, as we will come to see, the choice of boundary conditions do not have an affect on performance.
        
        We look for basis functions as a compact combination of Legendre polynomials,
        \be \label{eq:LG_general_basis}
        \phi_k(x)=L_k(x) + a_k L_{k+1} + b_k L_{k+2}(x),
        \ee
        where $L_k$ is the $k$-th Legendre polynomial and the parameters ${a_k,b_k}$ are chosen to satisfy the boundary conditions of the differential equations; see e.g. \cite{STW11} for more details.
        Such basis functions are referred to as modal basis functions.
        Hence, for the homogeneous Dirichlet boundary conditions we have $a_k=0$ and $b_k=-1$ in \eqref{eq:global_basis_func} which yields
        \be
            \phi_k(x)=L_k(x)-L_{k+2}(x).
            \label{eq:dirichlet_basis}
        \ee
        
        The set of global basis functions, $\{\phi_j\}$, are then used as test functions in the weak formulation of the derivative. 
        Noting that we have the zero boundary conditions at $x= \pm 1$, we multiply the DE in \eqref{eq:paradigm} by $\phi_j$, and integrate from $-1$ to $1$:
        \be \label{e:error1}
        	-\e \int^1_{-1}  {u}_{xx} \phi_j -  \int^1_{-1} {u}_x \phi_j \simeq  \int^1_{-1} f \phi_j.
        \ee
        By integration by parts, we find that
        \be
        	\e \int^1_{-1}  {u}_{x} (\phi_j)_x -  \int^1_{-1} {u}_x \phi_j \simeq \int^1_{-1} f \phi_j.
        \ee
        In this way, we can avoid the 2$^{nd}$ order derivative. 
        Hence, we set the residual
        \begin{gather}
            \begin{split}
                LHS := \sum_{j=0}^{m} \left( \e \int^1_{-1}  \widehat{u}_{x} {(\phi_j)}_x -  \int^1_{-1} \widehat{u}_x \phi_j\right),\qquad
                RHS := \sum_{j=0}^{m} \left( \int^1_{-1} f \phi_j \right),\\
                \mathcal{J}_{wf} := LHS - RHS, \hspace{4.1cm}
            \end{split}
        \end{gather}
        where $\widehat{u}$ is a numerical approximation of $u$,
        then minimize $\mathcal{J}_{wf}$ together with the other errors. 
        Here $m$ stands for the number of test functions of the weak formulation, which will be used in the loss function. 
        This loss function acts as a regularization factor \cite{bishop_2006, Tibshirani96}.
        For numerical derivatives and integration, we adopt spectral differentiations and Gauss-type integration formulas, which provide exponential convergence rate.
        The Gauss-type quadrature formulas provide powerful tools for evaluating not only integrals but also spectral differentiations \cite{RW05, DR07}.
        In addition, having precomputed the first-order differentiation matrix, the differentiation in the physical space can be carried out through matrix–matrix and matrix–vector multiplications.
        In fact, a similar approach was introduced by the authors in \cite{Raissi1026}. However, they used the weak form only for the loss function with an automatic derivative. 
        Without knowing this, we have developed our method based on the spectral methods. 
        In our method, we predict the coefficient of spectral approximation. When generating the training data set, these coefficients can be derived from Galerkin methods which are based on the weak formulation. Hence, it is natural to utilize the weak form to improve the convergence.

        We measure the accuracy of our predicted solution, $\widehat{u}$, against the ground truth solution, $u$, using a loss function, $\mathcal{L}$, defined by
        \be
            \mathcal{L}:= \mathcal{L}_{u} + \mathcal{L}_{wf},
        \ee
        % The measure of loss could be the mean absolute error, $MAE$, or the mean square error, $MSE$.
        where the loss corresponding to the numerical solution and weak form are $\mathcal{L}_u$ and $\mathcal{L}_{wf}$, respectively, and defined by
        \begin{align}
            \begin{split}
            \mathcal{L}_u&:=MSE\left\{u,\widehat{u}\right\},\\
            \mathcal{L}_{wf}&:=MSE\{LHS,RHS\},
            \end{split}
        \end{align}\\

       \subsubsection{Experimental Results}~\\
        
        \begin{figure}
            \begin{subfigure}[t]{.45\textwidth}
                \begin{center}
                    \includegraphics[width=0.9\textwidth]{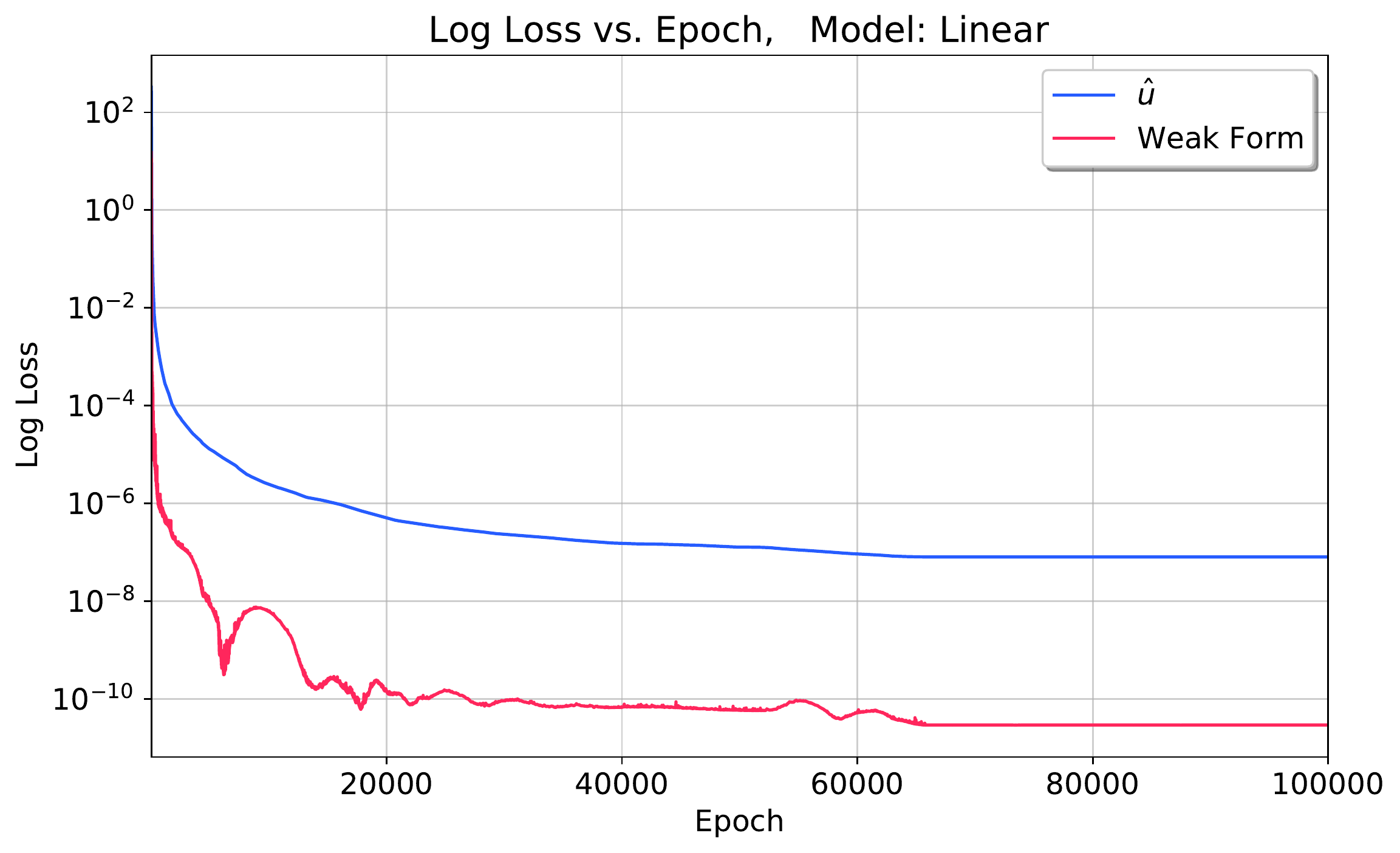}
                    \caption{Individual Losses}
                    \label{subfig:paradigm_loss_individual}
                \end{center}
            \end{subfigure}
            \hspace{0.5cm}
            \begin{subfigure}[t]{.45\textwidth}
                \begin{center}
                    \includegraphics[width=0.9\textwidth]{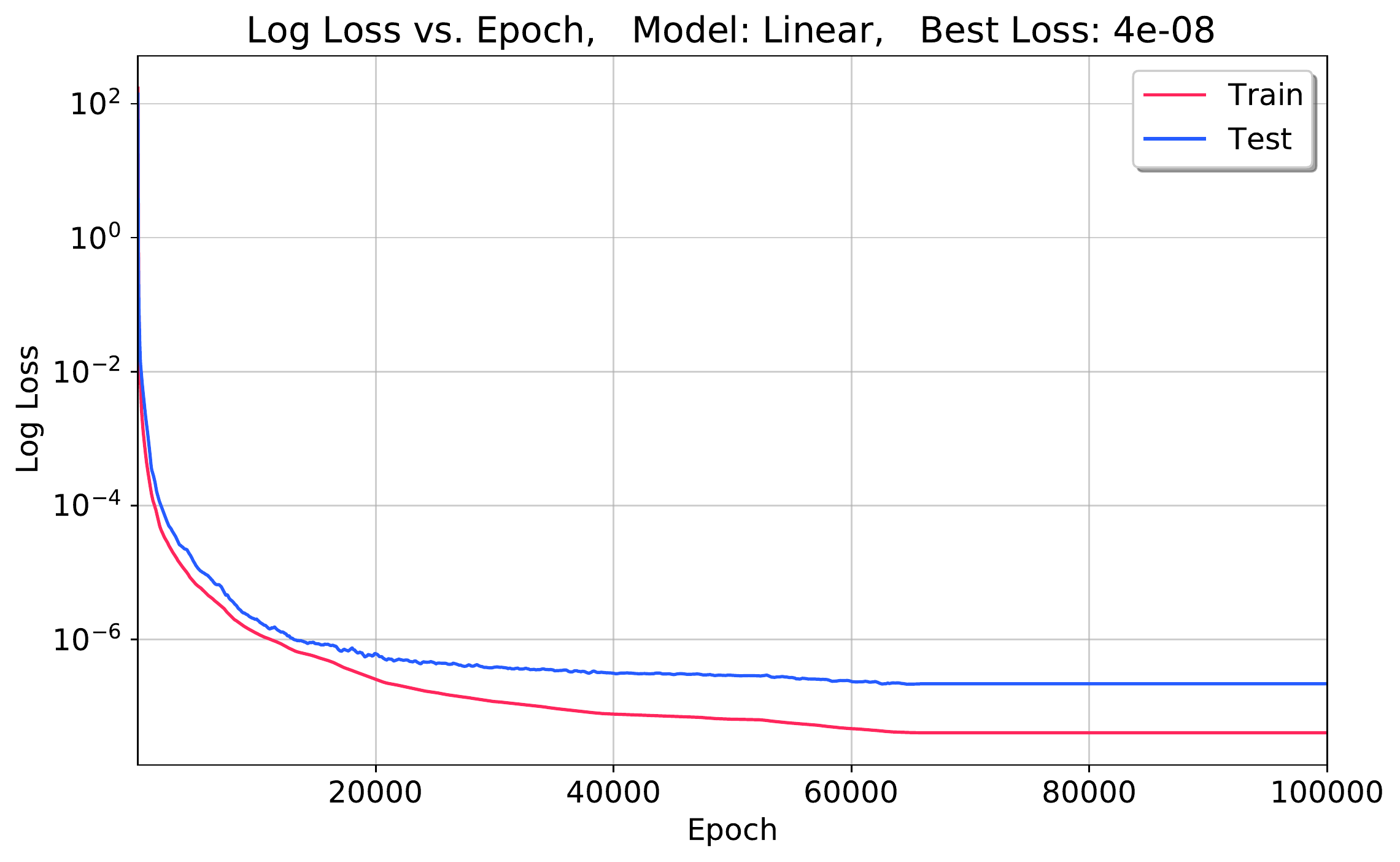}
                    \caption{Training Loss}
                    \label{subfig:paradigm_loss_train}
                \end{center}
            \end{subfigure}\\
            \begin{subfigure}[t]{.45\textwidth}
                \begin{center}
                    \includegraphics[width=0.9\textwidth]{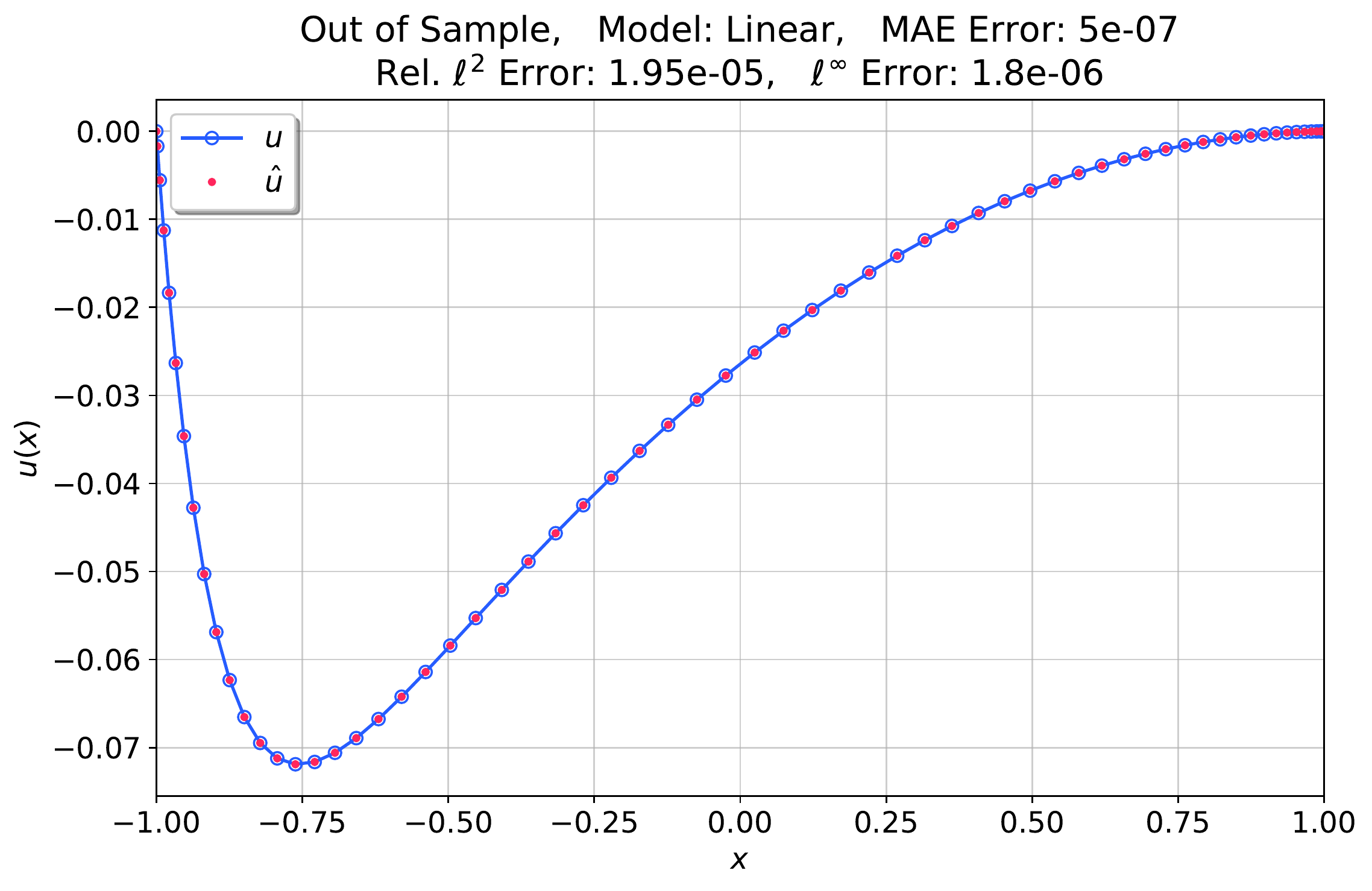}
                    \caption{Predicted Solution $\widehat{u}$}
                    \label{subfig:paradigm_sol_u}
                \end{center}
            \end{subfigure}
            \hspace{0.5cm}
            \begin{subfigure}[t]{.45\textwidth}
                \begin{center}
                    \includegraphics[width=0.9\textwidth]{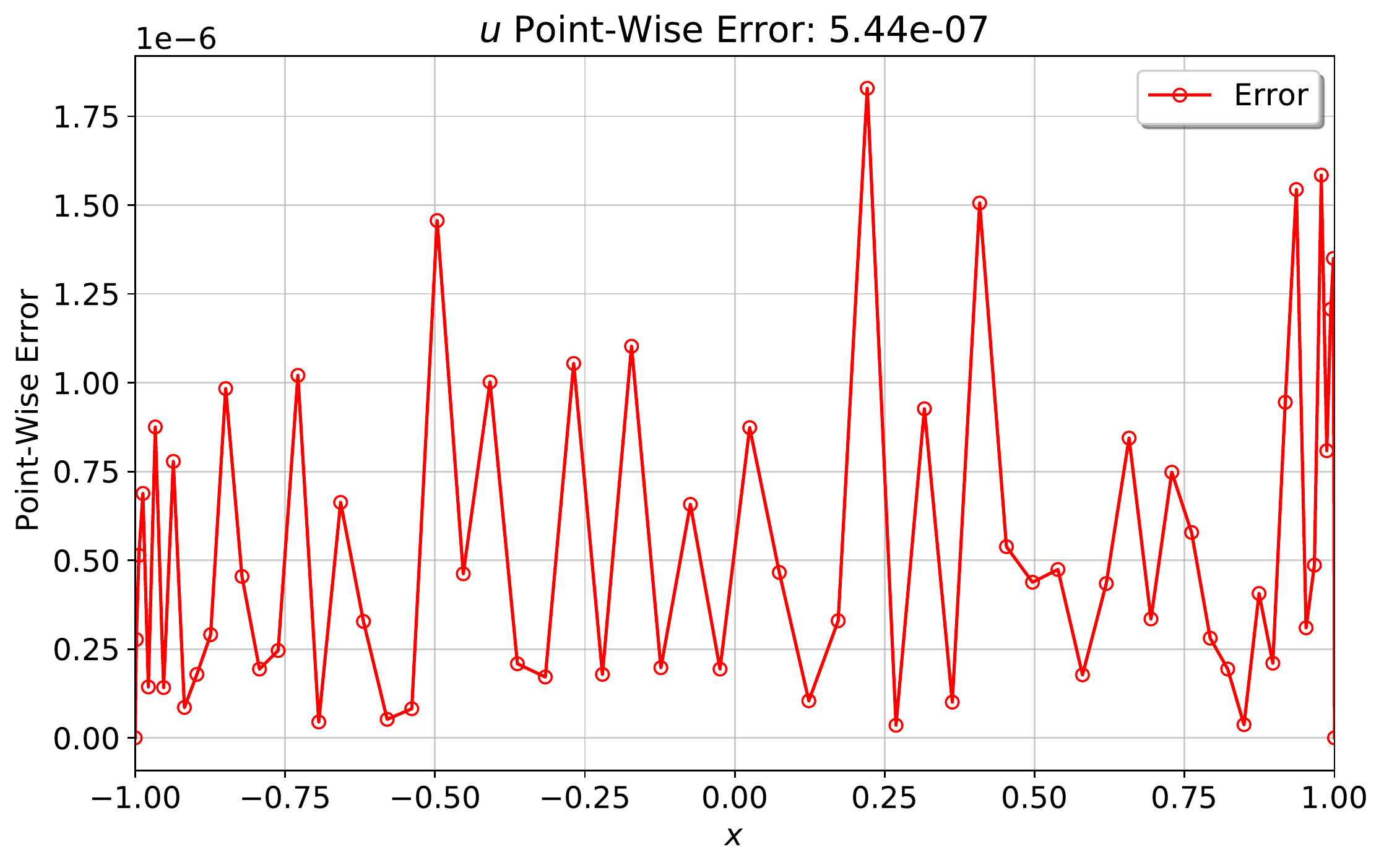}
                    \caption{Point-wise Error}
                    \label{subfig:paradigm_sol_upwe}
                \end{center}
            \end{subfigure}
            \caption{We trained a linear model \eqref{eq:paradigm} with $\e = 10^{-1}$ using a dataset with 10,000 solutions and 64 collocation points. The individual losses for $\mathcal{L}_u$ and $\mathcal{L}_{wf}$ are plotted on a semi-log plot over 100,000 epochs. The overall losses for the training and test sets can be seen on a semi-log plot beginning to separate as the number of epochs increases. An example of a predicted solution on out-of-sample data is observed with a mean relative $\ell^2$ error on the order of $10^{-5}$. The point-wise error plot shows a mean absolute error on the order of $10^{-7}$.}
            \label{fig:paradigm_sol}
        \end{figure} 

        Equation \eqref{eq:paradigm} can be solved with a {\it linear} architecture because it is a linear DE; the NN architecture can be seen in Figure \ref{subfig:architecture_linear}. Our architecture can achieve an {\it average} relative $\ell^2$ error on of order $10^{-5}$ on out of sample solutions with an extremely small data set; see Figure \ref{fig:paradigm_sol} for reference. For instance, by using a training data set with 100 randomly sampled solutions, our architecture is able to effectively predict sufficiently accurate solutions to Equation \eqref{eq:paradigm}. This actually speaks to the simplicity of a linear based differential equation and its representation; solutions to linear DEs can be represented better with linear networks than nonlinear models.
        
        For a linear model with 32 convolutional filters with a kernel size of 5, we were able to obtain an average relative $\ell^2$ error over an out-of-sample data set with 1,000 solutions on the order of $10^{-6}$. An example of an out-of-sample predicted solution can be seen in Figure \ref{subfig:paradigm_sol_u}. Our model generalizes well to out-of-sample data and is reflected in a $\log$-loss vs. epoch training loss plot in Figure \ref{subfig:paradigm_loss_train}. The individual losses were tracked in Figure \ref{subfig:paradigm_loss_individual} as a function of epoch and the loss associated with $\widehat{u}$ is a couple orders of magnitude greater than the loss associated with the weak form. The out-of-sample generalization is improved by increasing the size of input training set but the improvements are marginal at best.\\
        
        \subsubsection{Discussion}~\\
        
        Whether the training data set is $1,000$ solutions or $10,000$ solutions, the best performing models were on the order of $10^{-5}$ or $10^{-6}$, which suggests that there might be a fundamental limit in performance. This could be explored further in future works. By varying the architectural parameters the performance might increase marginally but work should be done to address this theoretically. Introducing a pooling operation could improve overall accuracy but we have not implemented that here. The performance could also improve by increasing the size of the training set but there is a trade-off being made because of the GPU memory limitations for training NNs one might have to switch from L-BFGS to a stochastic gradient method.
        
        To offer contrast, we did compare the performance of linear models to nonlinear models. Linear models converge rather quickly due to the L-BFGS optimization method but nonlinear models did not. Perhaps given enough time the models would converge but linear models outperformed nonlinear models, \textit{i.e.}, models with nonlinear activation functions do not perform well on linear DEs. All nonlinear models that we trained and evaluated drastically under-performed linear models \cite{lecun_gradient}.\\

    \subsection{Helmholtz Equation}~\\
    
        Next, we will utilize our algorithm to find solutions to the Helmholtz equation with Neumann boundary conditions. This equation offers a contrast to \eqref{eq:paradigm} because it is a form of the wave equation with varying use in the field of optics \cite{steward2004fourier}. Also, different boundary conditions are used, furthering the notion of the general applicability of our algorithm. 
        Consider the Helmholtz differential equation given as
        \be
            \label{eq:helmholtz}
            \begin{cases}
            	 u_{xx} + k_u u = f(x),\\
            	u'(- 1) = u'(1) = 0,
            \end{cases}
        \ee
        where $f$ is defined as in \eqref{eq:forcing}.
        We remark that the governing equation, \eqref{eq:helmholtz}, describes essential aspects of scattering of linear waves by periodic multiply layered gratings. Indeed, considering the quasiperiodicity of solutions, the generalized Fourier (Floquet) series expansions reduce the transverse electric and magnetic waves to the one-dimensional inhomogeneous Helmholtz equations, which can be solved by numerical approaches; see e.g. \cite{HN17, HN18}.
        
        The global basis function (Legendre basis function) of \eqref{eq:helmholtz} is defined by the homogeneous Neumann boundary conditions. Hence, in \eqref{eq:LG_general_basis}, we set 
        $$a_k=0,\quad\text{and}\quad b_k= -\dfrac{k(k+1)}{(k+2)(k+3)}$$
        in \eqref{eq:global_basis_func} which yields
        \be
            \phi_k(x)=L_k(x) -\dfrac{k(k+1)}{(k+2)(k+3)} L_{k+2}(x),
        \ee
        where $0\leq k \leq N-2$; for more details, see e.g. \cite{STW11}.
        
        The set of global basis functions, $\{\phi_k\}$, are then used as test functions in the weak formulation,
        \be
        - \int^1_{-1} u_x (\phi_k)_x dx + k_u \int^1_{-1} u \phi_k dx
        = \int^1_{-1} f \phi_k dx.
        \ee
        Hence, we define
        \begin{gather}
            \begin{split}
                LHS := \sum_{j=0}^{m} \left( - \int^1_{-1}  \widehat{u}_{x} {(\phi_j)}_x +  \int^1_{-1} k_u \widehat{u} \phi_j\right), \qquad 
                RHS  := \sum_{j=0}^{m} \left( \int^1_{-1} f \phi_j \right),\\
                \mathcal{J}_{wf} := LHS - RHS,\hspace{4cm}
            \end{split}
        \end{gather}
        then minimize $\mathcal{J}_{wf}$ together with the other errors. 
                Here $m$ stands for the number of test functions of the weak formulation, which will be used in the loss function. 
        The loss function, $\mathcal{L}$, for our optimal model was
        \be
            \mathcal{L}=\mathcal{L}_u+\mathcal{L}_{wf},
        \ee
        where $\mathcal{L}_u$ and $\mathcal{L}_{wf}$ as defined in \eqref{eq:mse}.\\

    \subsubsection{Experimental Results}~\\
    
        \begin{figure}
            \begin{subfigure}[t]{.45\textwidth}
                \begin{center}
                    \includegraphics[width=0.9\textwidth]{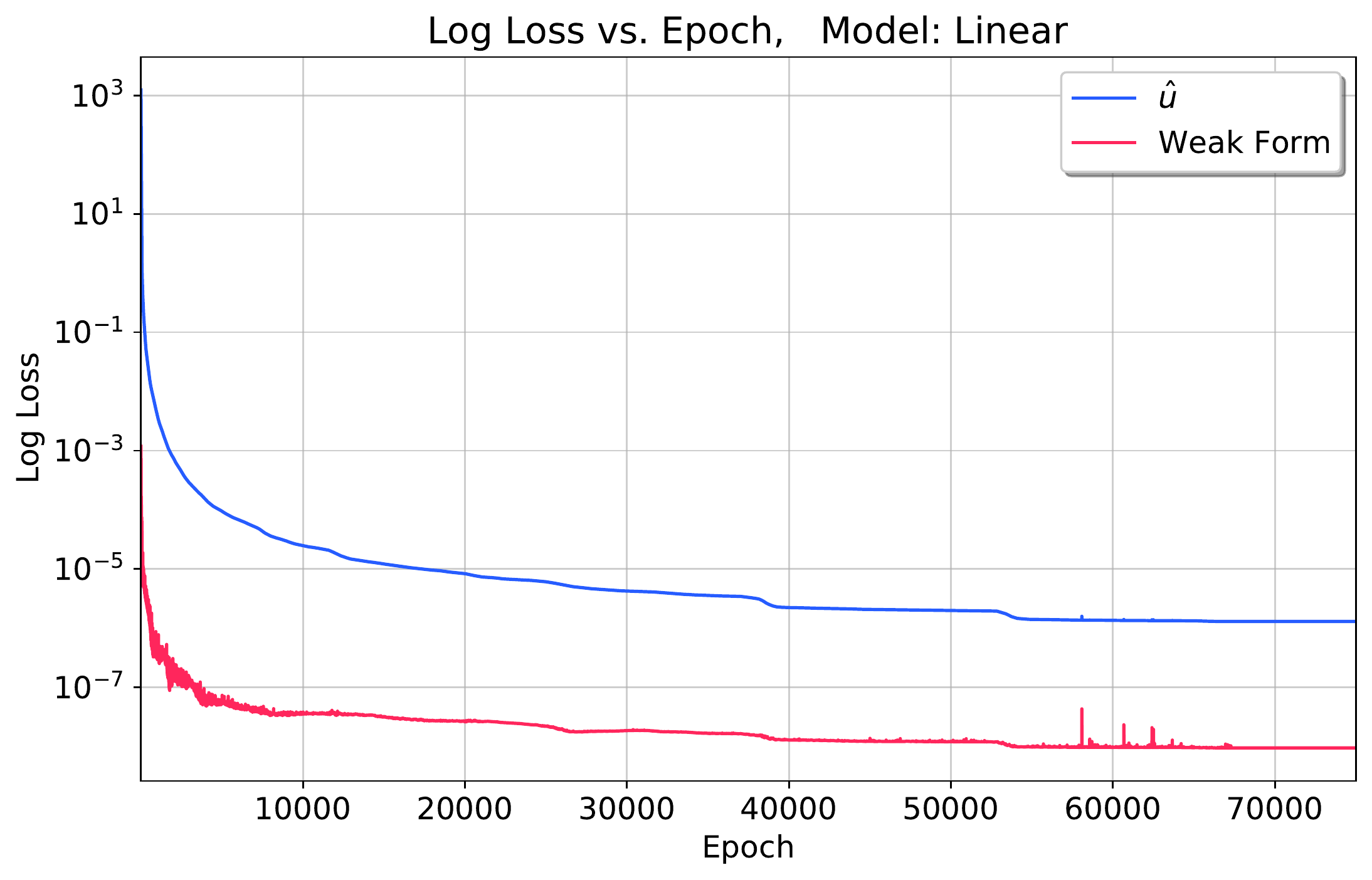}
                    \caption{Individual Losses}
                    \label{subfig:helmholtz_loss_individual}
                \end{center}
            \end{subfigure}
            \hspace{0.5cm}
            \begin{subfigure}[t]{.45\textwidth}
                \begin{center}
                    \includegraphics[width=0.9\textwidth]{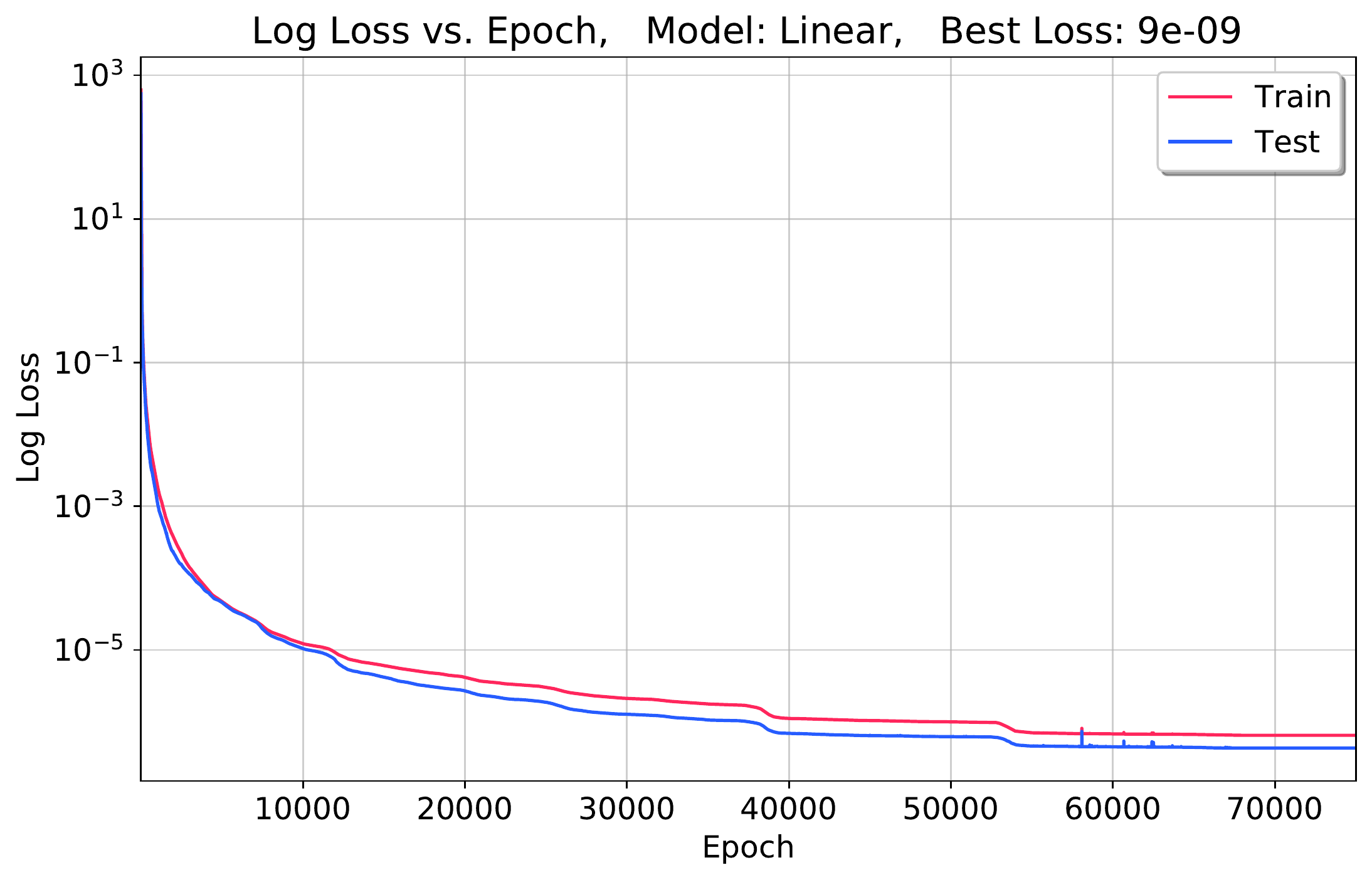}
                    \caption{Training Loss}
                    \label{subfig:helmholtz_loss_train}
                \end{center}
            \end{subfigure}\\
            \begin{subfigure}[t]{.45\textwidth}
                \begin{center}
                    \includegraphics[width=0.9\textwidth]{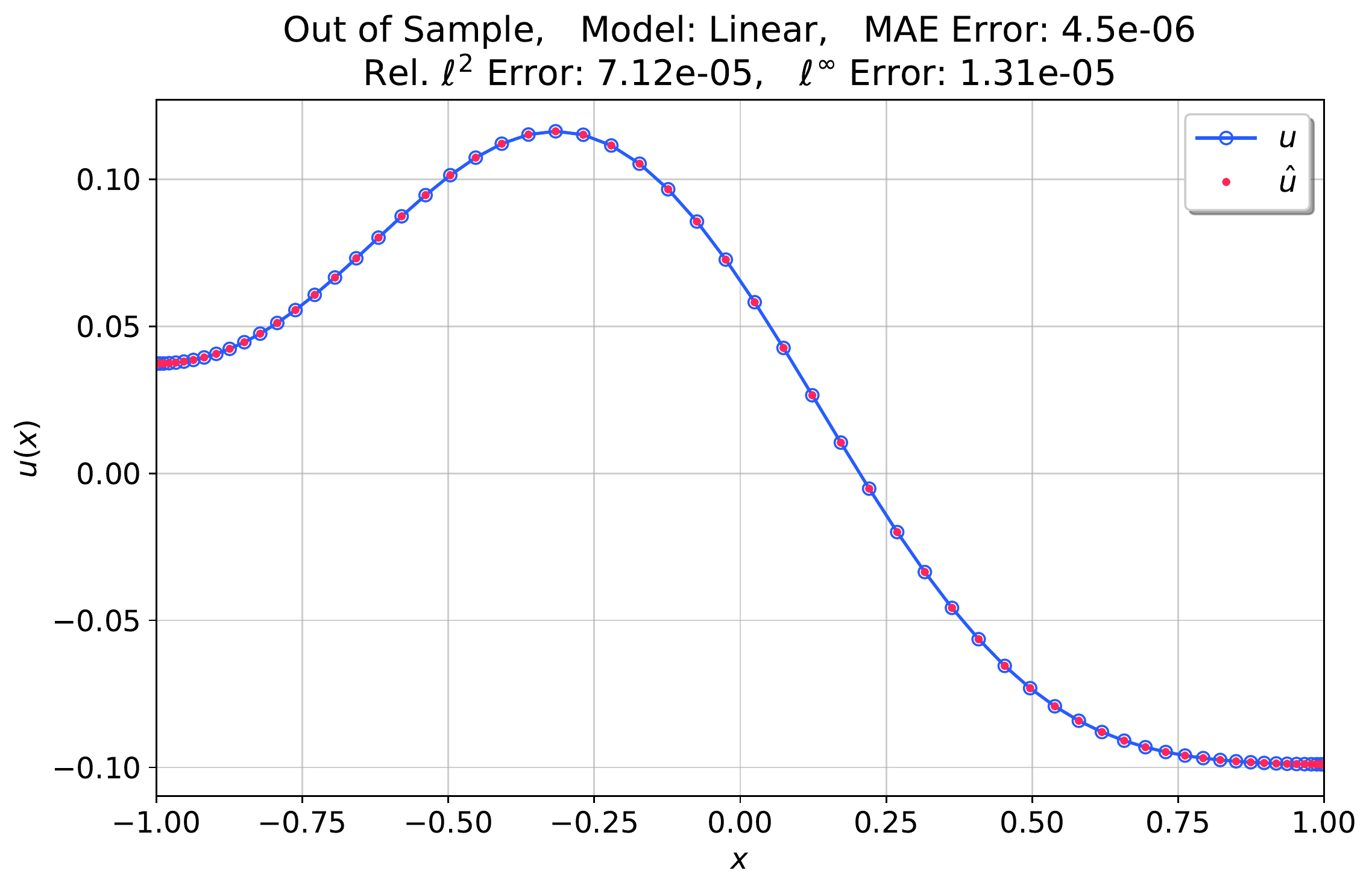}
                    \caption{Predicted  solution $\widehat{u}$}
                    \label{subfig:helmholtz_sol_u}
                \end{center}
            \end{subfigure}
            \hspace{0.5cm}
            \begin{subfigure}[t]{.45\textwidth}
                \begin{center}
                    \includegraphics[width=0.9\textwidth]{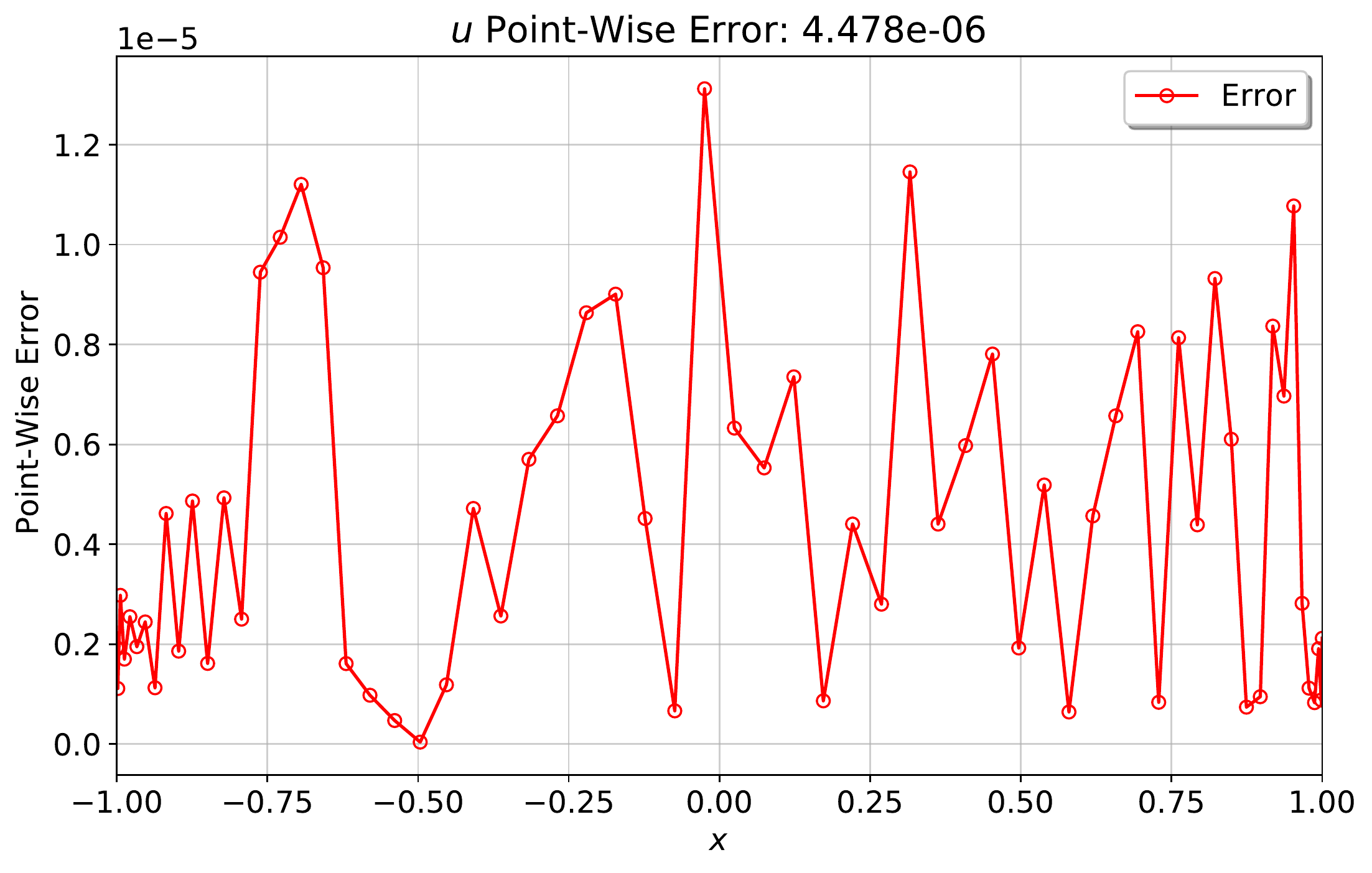}
                    \caption{Point-wise Error}
                    \label{subfig:helmholtz_sol_upwe}
                \end{center}
            \end{subfigure}
            \caption{Here is a trained model \eqref{eq:helmholtz} with $k_u = 3.5$ using 10,000 solutions in the training set with 64 collocation points. The individual losses for $\mathcal{L}_{u}$ and $\mathcal{L}_{wf}$ are plotted on a semi-log plot over 100,000 epochs, but this model converged after approximately 65,000 epochs. The overall losses for the training and test sets can be seen on a semi-log plot, where we observe the test loss is \textit{less than} the training loss. An out-of-sample predicted solution is plotted along with the known true solution with a mean relative $\ell^2$ error on the order of $10^{-5}$. A plot of the point-wise error is observed where the mean absolute error is on the order of $10^{-5}$.}
            \label{fig:helmholtz_sol}
        \end{figure}
        
        The Helmholtz equation is a linear wave equation with the optimal neural architecture again being entirely linear. We use the same architecture as seen in Figure \ref{subfig:architecture_linear} to achieve similar results to the CDE. A purely linear neural architecture, without any nonlinear activation functions, is able to generalize to a larger test set from a small training set. This is achieved by training the model for `many' epochs until either the model converges, or resources need to be freed for other projects. A `converged' model in this context is one that is defined to be where the loss function, $\mathcal{L}$, has converged, \textit{i.e.}, the loss at some epoch is the same as the next epoch, \textit{ad infinitum}.
        
        A linear model with 10,000 in-sample solutions and 64 collocation points was trained for 100,000 epochs, where the results can be seen in Figure \ref{fig:helmholtz_sol}. The loss values for the predicted solution, $\widehat{u}$, and its associated weak form were tracked at each epoch and can be seen to converge in Figure \ref{subfig:helmholtz_loss_individual}. Testing was performed with a 1,000 out-of-sample data set, as seen in Figure \ref{subfig:helmholtz_loss_train}, where we observe that both the training and test sets smoothly converge together indicating that our model generalizes rather well to out-of-sample data for small input training set. Figure \ref{subfig:helmholtz_sol_u} shows a predicted solution, $\widehat{u}$, compared to the true solution generated with the SEM with a relative $\ell^2$-error of 7.12$\times 10^{-5}$.\\
        
    \subsubsection{Discussion}~\\
    
        Linear equations can be solved with a linear architecture but nonlinear equations cannot be solved with linear architectures. This has multiple benefits, such as reduction in computational complexity, the ability to generalize from a small training data set to a larger out of sample data set, \textit{etc}, but it has its limitations. All linear models, as depicted in Figure  \ref{subfig:architecture_linear}, were created with a kernel size of 5, stride of 1, padding 2 and 32 filters per channel, similar to Bar-Sinai et al 2019 \cite{Bar_Sinai_2019}. The LGNet algorithm finds accurate solutions to Equations \eqref{eq:paradigm} and \eqref{eq:helmholtz} regardless of the type of boundary conditions. We also considered different hyperparameters such as varying the kernel size, number of convolutional filters per layer, but were unable to find a discernible pattern.
        
        It is important to note the difference between the point-wise error plots in Figures \ref{subfig:paradigm_sol_upwe} and \ref{subfig:helmholtz_sol_upwe}. In Equation \eqref{eq:paradigm} we use homogeneous Dirichlet boundary conditions and in Equation \eqref{eq:helmholtz} we use homogeneous Neumann, which is reflected in the numerical error difference plots. In Figure \ref{subfig:paradigm_sol_upwe} the model is able to properly learn the boundary conditions to a high degree of accuracy but the boundary condition generates numerical errors at the boundary in Figure \ref{subfig:helmholtz_sol_upwe}. Although the average relative $\ell^2$-norm error over the test set is comparable, there is a distinguishing difference between models with Dirichlet and Neumann boundary conditions, but the accuracy of both models are comparable. The performance of each model is, roughly speaking, the same regardless of which boundary conditions are used.

    \section{Non-linear model} \label{s:4}
    The performance of neural networks, when applied to linear differential equations, prove to be powerful and accurate.
    As a natural extension, we consider a neural network for non-linear models.
    Nonlinear differential equations are notoriously difficult to solve and we should not assume that a linear model would work for a nonlinear differential equation. Our neural architecture should, in a sense, reflect the differential equation, \textit{i.e.}, we want to use an architecture with nonlinear operations to represent the solution of a nonlinear equation. In this section we will introduce nonlinear neural architectures with various activation functions to explore the accuracy of trained models.\\

    \subsection{Burgers Equation}~\\
    
        This problem can be further investigated by applying our algorithm to the canonical problem known as Burgers equation; a differential equation known for its nonlinear shockwave characteristics. We define Burgers equation as
        \be
            \label{eq:Burgers}
            \begin{cases}
            	-\e u_{xx} + u u_x = {f},\\
            	g(\pm 1) = 0,
            \end{cases}
        \ee
        where $\e$ is a diffusion coefficient. 
        %In the previous models we used an input training set the was generated using with normally distributed parameters and obtained a nice accuracy. The same does not occur for Equation \eqref{eq:Burgers}. In the new modified forcing function we switch to a uniform distribution, where $\tilde{f}$ is defined by 
        To avoid a very small solution, which is not interesting, we utilize a uniform distribution when we generate random coefficients on the external forcing function ${f}$, that is our input data:
        \begin{equation}
            {f}=(3+\upsilon_1)\sin[(1+\upsilon_2)\pi x]+(3+\upsilon_3)\cos[(1+\upsilon_4)\pi x],
            \label{eq:tilde_f}
        \end{equation}
        with $\upsilon_i\in {\text{Uniform}}[0,2]$ for $1 \leq i\leq 4$.
        
        As an input training data set with a large variance performed poorly, we correct this by normalizing our input training set to mean 0 and standard deviation 1, similar data augmentation was utilized in \cite{Bar_Sinai_2019}. Data pre-processing is a common technique used in machine learning for improved network performance. Please see \cite{bishop_2006} or \cite{vapnik_2000} for more on data pre-processing.
        
        Since homogeneous Dirichlet boundary conditions are used, the global basis functions for \eqref{eq:Burgers} are given by \eqref{eq:dirichlet_basis}. 
        We remark that the proposed algorithm can compute other boundary conditions, such as Neumann or Robin boundary conditions, since the Legendre basis gives an exact boundary condition.
        We can similarly derive the weak form as laid out above by multiplying each side of \eqref{eq:Burgers} by some test function, $\phi_k$, and integrating from -1 to 1 to yield
        \be
        \begin{split}
            & \,\, -\e \int^1_{-1} u (\phi_k)_{xx}\,dx + \frac{1}{2}\int^1_{-1}  (u^2)_x (\phi_k)\, dx = \int^1_{-1} f \phi_k\, dx \\
            & \Rightarrow
            \e \int^1_{-1} u_x (\phi_k)_x\,dx - \frac{1}{2}\int^1_{-1}  u^2 (\phi_k)_x\, dx = \int^1_{-1} f \phi_k\, dx.
        \end{split}
        \ee
        We again set the residual as
        \begin{gather}
            \begin{split}
                LHS := \sum_{j=0}^{m} \left( \e \int^1_{-1} u_x (\phi_j)_x\,dx - \frac{1}{2}\int^1_{-1}  u^2 (\phi_j)_x\, dx\right),\quad
                RHS := \sum_{j=0}^{m} \left( \int^1_{-1} f \phi_j\, dx \right),\\
                \mathcal{J}_{wf} := LHS - RHS.\hspace{6cm}
            \end{split}
        \end{gather}
        Being nonlinear, Burgers' equation poses an interesting problem for the utilization of NNs to solve differential equations. In order to create a training data set we had to utilize the Picard iteration method with a tolerance level of $10^{-9}$. The limitations of a nonlinear data set are met with less stringent requirements with regards to training set data, but we should note that even the best linear models could not achieve such accuracy, so this tolerance will be appropriate for training set purposes.\\
    
    \subsection{Experimental Results}~\\
        \begin{figure}[t]
            \begin{center}
                \includegraphics[width=0.9\textwidth]{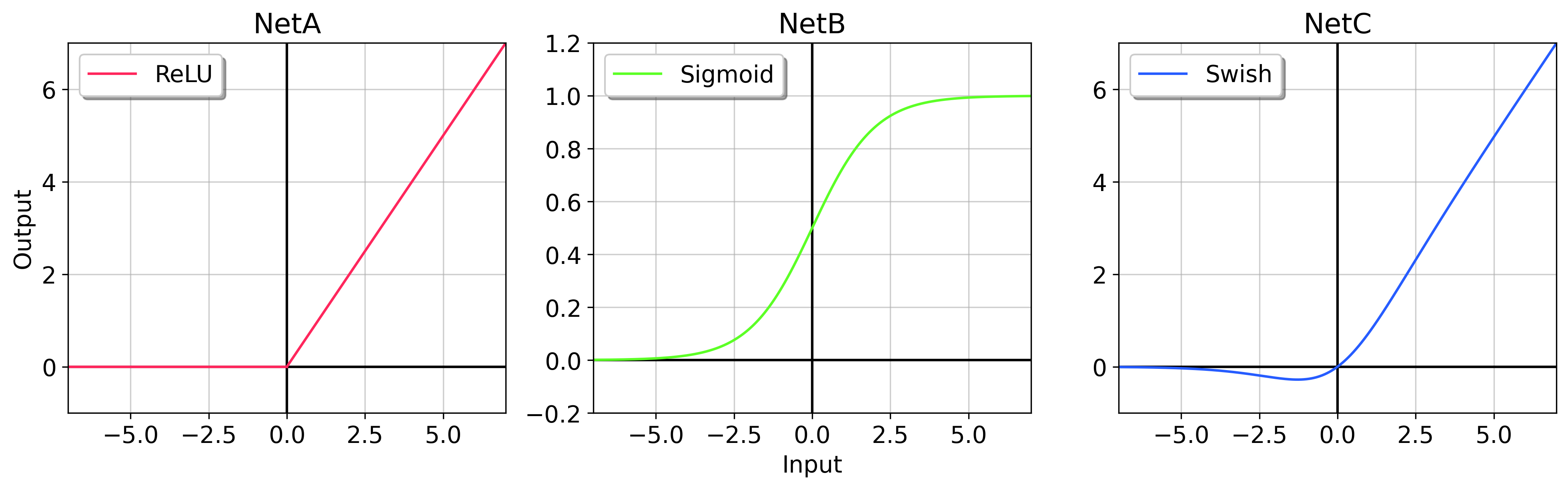}\\
                \includegraphics[width=0.9\textwidth]{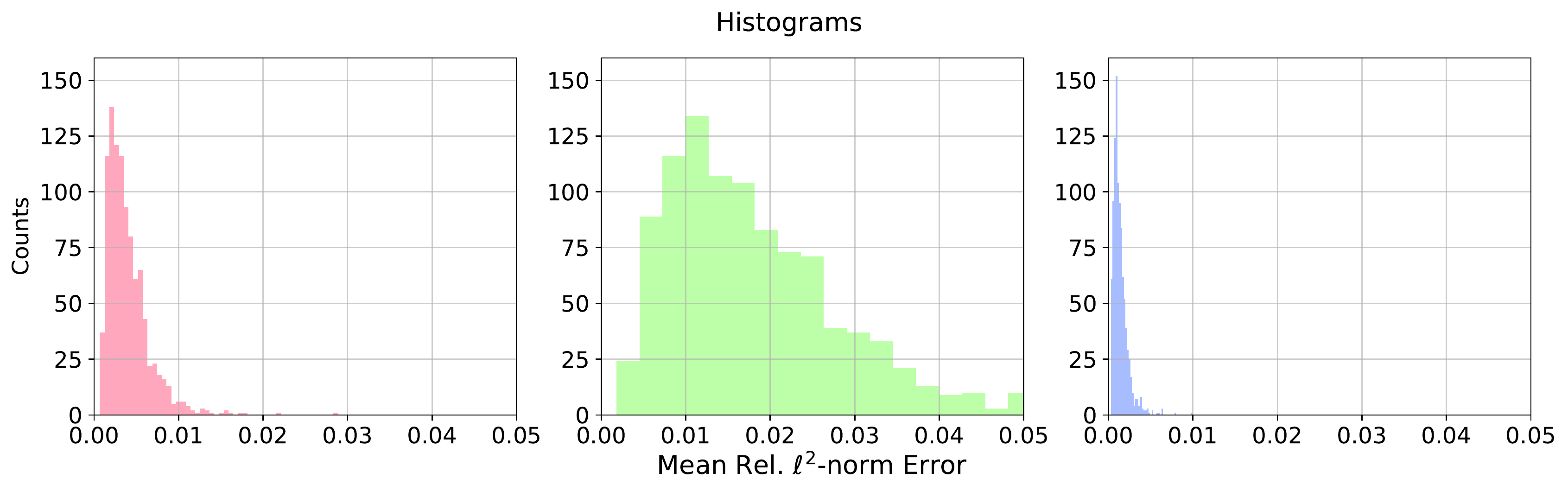}\\
                \caption{We investigate the effect that different activation functions have on the performance of LGNet. (Top) The activation function for each network is plotted. (Bottom) The relative $\ell^2$-norm error histogram for each sample in the test dataset is recorded and a histogram plot is shown below its activation function. Each model is trained with an input training set of 15,000 solutions, 64 collocation points, 32 filters, 5 kernel, and 5 convolutional blocks. Models trained with the \textit{Swish} activation function consistently outperform models using ReLU or sigmoid.}
                \label{fig:activations}
            \end{center}
        \end{figure}
        
        \begin{figure}
            \begin{center}
                \includegraphics[scale=0.375]{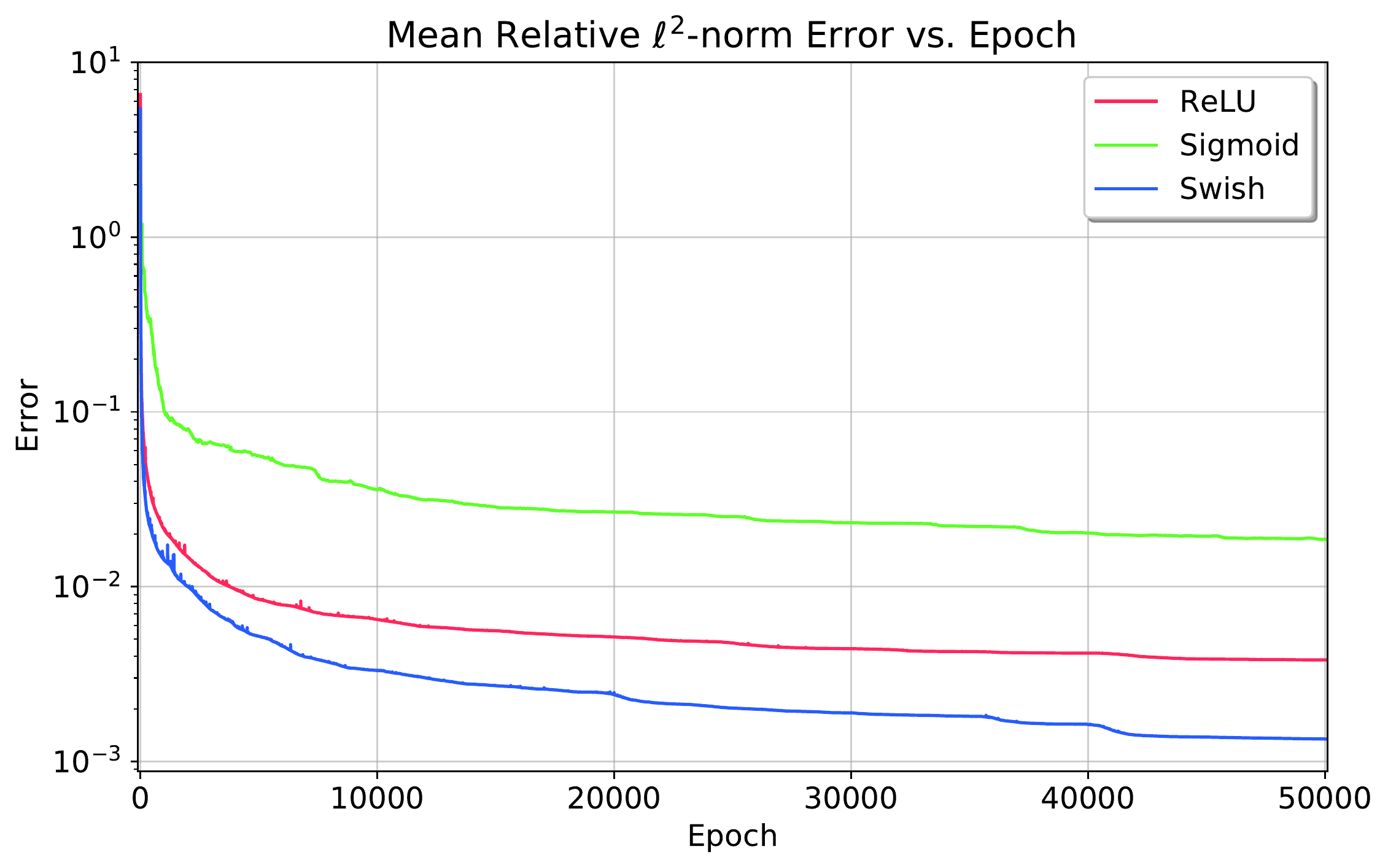}
                \begin{minipage}{0.8\linewidth}
                	\caption{A comparison of different activation functions can be seen in the mean relative $\ell^2$-norm error vs. epoch plot above; each model has the same parameters as Figure \ref{fig:activations}. The error is computed on the test data set and recorded at each epoch. The \textit{Swish} activation function function consistently outperforms both sigmoid and ReLU. Both ReLU and \textit{Swish} are an order of magnitude better than sigmoid, on average. The performance of each model is a reflection of the histograms in Figure \ref{fig:activations}.}
                	\label{fig:errors}
                \end{minipage}
            \end{center}
        \end{figure}
          
        \begin{figure}
            \begin{subfigure}[t]{.45\textwidth}
                \begin{center}
                    \includegraphics[width=0.9\textwidth]{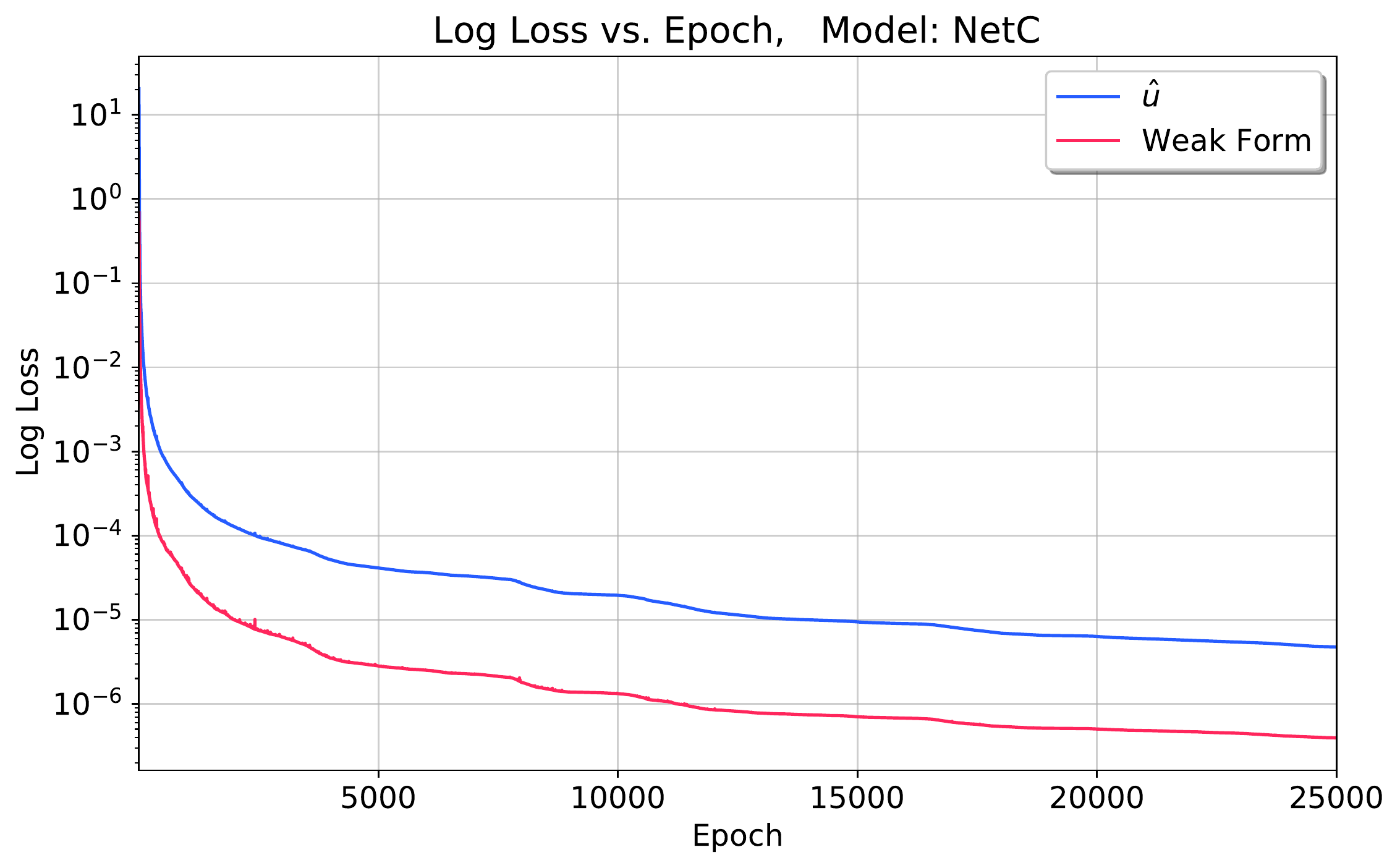}
                    \caption{Individual Losses}
                    \label{subfig:burgers_loss_individual}
                \end{center}
            \end{subfigure}
            \hspace{0.5cm}
            \begin{subfigure}[t]{.45\textwidth}
                \begin{center}
                    \includegraphics[width=0.9\textwidth]{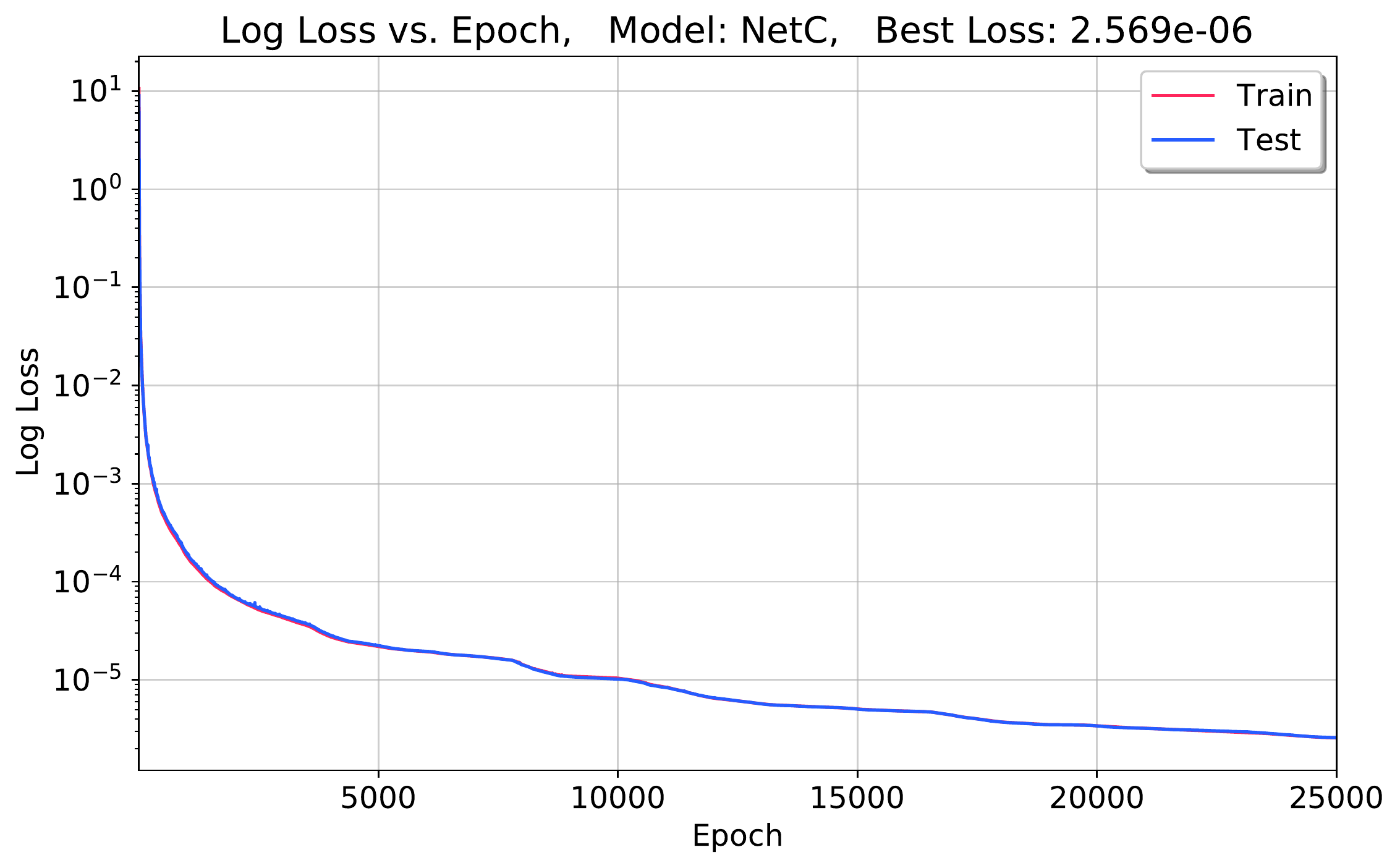}
                    \caption{Training loss}
             			\label{subfig:burgers_loss_train}
                \end{center}
            \end{subfigure}\\
            \begin{subfigure}[t]{.45\textwidth}
                \begin{center}
                    \includegraphics[width=0.9\textwidth]{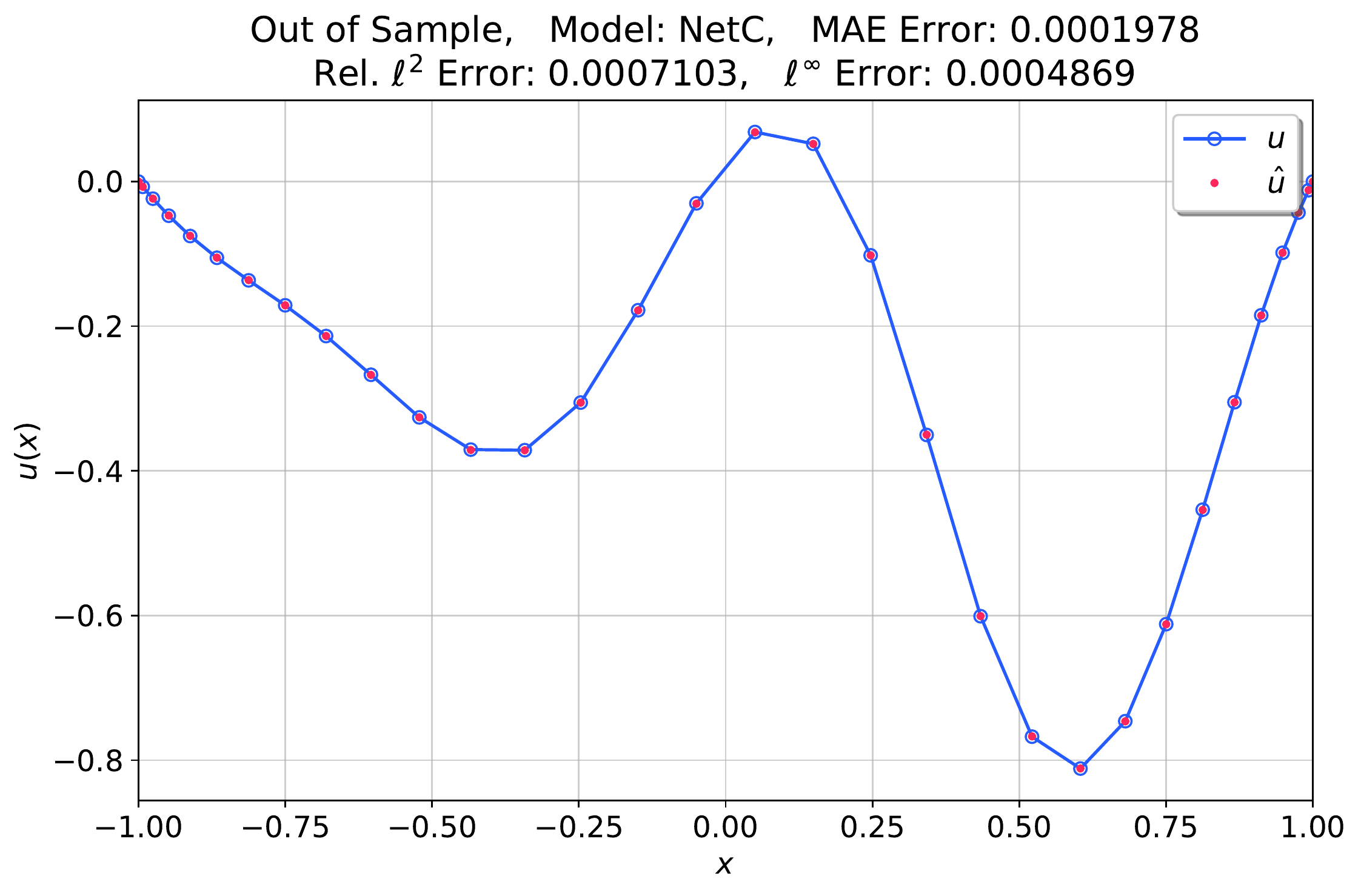}
                    \caption{Predicted solution $\widehat{u}$}
                    \label{subfig:burgers_sol_u}
                \end{center}
            \end{subfigure}
            \hspace{0.5cm}
            \begin{subfigure}[t]{.45\textwidth}
                \begin{center}
                    \includegraphics[width=0.9\textwidth]{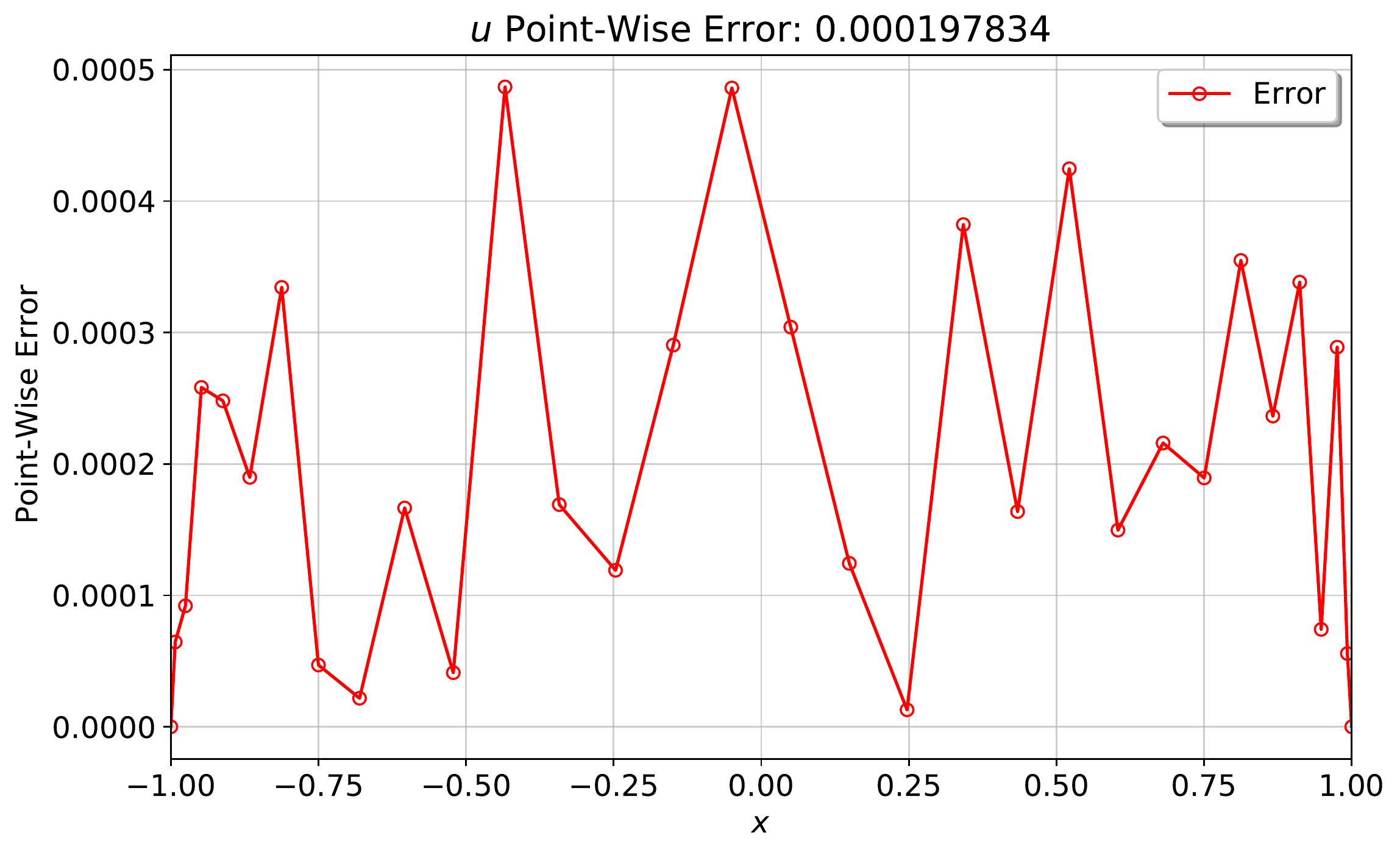}
                    \caption{Point-wise Error}
                    \label{subfig:burgers_sol_upwe}
                \end{center}
            \end{subfigure}
            \caption{Here is a model of \eqref{eq:Burgers} with $\varepsilon = 0.5$ using 10,000 solutions in the training set with 31 collocation points. The individual losses for $\mathcal{L}_u$ and $\mathcal{L}_{wf}$ are plotted on a semi-log plot and the overall losses for the training and test sets can be seen on a semi-log plot. An out-of-sample predicted solution is plotted along with the known true solution with a mean relative $\ell^2$ error on the order of $10^{-4}$. A plot of the point-wise error is observed where the average mean absolute error is on the order of $10^{-4}$.}
            \label{fig:burgers_sol}
        \end{figure}

        Initially, we use a neural architecture similar to the one found in \cite{Bar_Sinai_2019, zhuang2020learned} to start with, as seen in Figure \ref{subfig:architecture_netA} and referred to as \texttt{NetA}. Current state-of-the-art NNs, such as deep residual networks, show an increase in performance as the depth of the network increases, but ResNet did not perform well when using either Group Norm, or Batch Norm \cite{he2015deep, Ioffe15, wu2018group}. We did observe that the $\texttt{NetA}$ style of architecture was optimal when compared to ResNet for the LGNet. Each convolutional has 32 filters that are learned with a kernel size of 5 and a stride of 1 and paddding of 2. The depth of this network is increased by adding additional blocks, $\mathbb{B}$.
        
        For reference, we also implement the sigmoid activation function as a drop-in replacement for ReLU by comparison; we call this neural architecture $\texttt{NetB}$ and its activation can be seen in Figure \ref{fig:activations}. ReLU was seen to consistently outperform sigmoid by an order of magnitude. Saturation effects have been previously mentioned with regards to sigmoid in literature but we simply note a poor performance with regards to ReLU.
        
        Recent work has demonstrated the effectiveness of the `\textit{Swish}' activation function; see Figure \ref{fig:activations} \cite{ramach2017searching}. This activation function has shown consistent improvement over the ReLU activation function in the application of CNNs to many different data sets which reflected in Figure \ref{fig:errors}. We use this as a drop-in replacement in our $\texttt{NetA}$ architecture, for comparison, and refer to this network as $\texttt{NetC}$. The \textit{Swish} activation function is given by
        \be
        g(x):=x\cdot \text{sigmoid}(x)=\frac{x}{1+e^{-x}},
        \label{eq:swish}
        \ee
        and does offer a consistent improvement over the ReLU activation function. Equation \eqref{eq:swish} is a non-monotonic one-sided unbounded smooth function, which will not experience saturation, a drawback of the sigmoid activation function in \texttt{NetB}. Furthermore, the \textit{Swish} activation function offers both negative and positive activation values, which is a characteristic that ReLU does not have.
        
        An example of a trained \texttt{NetC} model can be seen in Figure \ref{fig:burgers_sol}. The parameters of our best model were 32 filters, kernel size of 5, and 10 blocks. The model is cross-validated by comparing the trained model with out-of-sample data at each epoch as seen in Figure \ref{subfig:burgers_loss_train}. Again, we used the weak form for regularization and the loss values can be seen in Figure \ref{subfig:burgers_loss_individual}. A predicted solution, $\widehat{u}$, can be seen in Figure \ref{subfig:burgers_sol_u} for an out-of-sample forcing function, $\tilde{f}$, with a relative $\ell^2$-norm error of $\sim 7.1\times 10^{-4}$. The average point-wise error for this solution is $\sim 2\times 10^{-4}$ and can be seen in Figure \ref{subfig:burgers_sol_upwe}. The mean relative $\ell^2$-norm error on an out-of-sample data set of 1,000 solutions is $8.8\times 10^{-4}$.
        
        We have tabulated the numerical results for each equation in Table \ref{tab:results}. For CDE we used a linear model and let the model train for 100,000 epochs. The model converged at $\sim$ 70,000 epochs and achieved a mean relative $\ell^2$-norm error on the order of $10^{-6}$. Similarly, we used a linear architecture for the Helmholtz equation and achieved a mean relative $\ell^2$-norm error on the order of $10^{-5}$. Finally, our best model for Burgers' equation utilized the \texttt{NetC} architecture with 25,000 epochs. This model was able to achieve a mean relative $\ell^2$-norm error on the order of $10^{-4}$.\\
    
        \begin{table}    
\centering
    \begin{tabular}{||c||c|c|c|c|c|c||c||}
        \hline
        Equation & Model & Epochs & blocks & filters & ks & Mean Rel. $\ell^2$ Error \\ \hline \hline
        CDE & Linear & 100,000 & - & 32 & 5 & 1.26$\times 10^{-6}$ \\ \hline
        Helmholtz & Linear & 100,000 & - & 32 & 5 & 4.57$\times 10^{-5}$ \\ \hline
        Burgers & \texttt{NetC} & 25,000 & 4 & 32 & 5 & 8.81$\times 10^{-4}$ \\ \hline
    \end{tabular}
    \caption{The best trained models are presented here. For each equation, Convection Diffusion, Helmholtz, and Burgers, the parameters used in training are tabulated for clarity. Each model can achieve a mean relative $\ell^2$-norm error of $\mathcal{O}(10^{-5})$ when using the L-BFGS optimizer. Data augmentation is used on nonlinear input data to realize comparable mean relative $\ell^2$-norm error. Thus, effectively gaining comparable performance on nonlinear equations as linear equations.}
     \label{tab:results}
\end{table}
    
    \subsection{Discussion}~\\
    
        Nonlinear equations are challenging with regards to numerical methods, as well as NNs, and this problem can be resolved with data augmentation. The statistical nature of the neural architectural model introduces more room for error when training NNs, if the variance of the input data set is too large, but still yields sufficiently accurate results. The L-BFGS algorithm has its benefits during full-batch gradient descent, but when the size of the data set exceeds the memory limitations of the GPU, mini-batch L-BFGS gradient descent can diverge; thus, switching to SGD, or Adam, as an optimizer might yield similar results \cite{kingma2017adam, ruder2016overview, saad_1998}. Also, introducing a pooling operation into the neural architecture could greatly improve upon local invariance, but we have not measure the effects here.
        
        The performance of a DNN to predict solutions of Equation \eqref{eq:Burgers} can yield comparable results to Equations \eqref{eq:paradigm} or \eqref{eq:helmholtz}, but only with data augmentation. This was used in Bar-Sinai et al. 2019, but we note that data augmentation on linear equations is not necessary. If the input forcing function, ${f}$, in Equation \eqref{eq:tilde_f} is normalized prior to training then nonlinear equations can enjoy an accuracy comparable to that of linear equations.
        
        \begin{rem}
            The current architecture can be extended to time dependent problems, a PDE model. For instance, we consider a time dependent Burgers' equations,
            \be \label{e:time}
            \begin{split}
                 u_t - u_{xx} + u u_x &= f, \quad x \in (-1,1),\\
                 u(\pm 1) &= 0.
            \end{split}
            \ee
            Applying the backward Euler method in time, the equation $\eqref{e:time}_1$ becomes
            \be
                \dfrac{u^{n+1} - u^{n}}{\Delta t}
                - u^{n+1}_{xx} + u^{n+1} u_x^{n+1} = f^{n+1}.
            \ee
            Hence, by setting $v:= u^{n+1}$ we deduce that
            \be \label{e:new_time}
                -v_{xx} + (v^2)_x + \f{1}{\Delta t} v = f^{n+1} + \f{1}{\Delta t} u^n.
            \ee
            However, the structure of \eqref{e:new_time} is very similar to that in Section \ref{s:4}, and a similar approach/architecture of Section 4 can be used.
            It is noteworthy to mention that we need to build only one network for all time steps. 
            Indeed, the only difference in different time steps is $u^n$ in \eqref{e:new_time}, which is our input data.
            Hence, one neural network is enough to solve the time-dependent model.
            Further discussion will be provided in the forthcoming article.
        \end{rem}

    \section{Conclusion} \label{s:5}
    Spectral element methods are used to create accurate data sets, then our LGNet algorithm is able to find suitably accurate solutions. We use deep neural networks to \textit{learn} the coefficients relative to some global basis then are able to accurately reconstruct the numerical solution. We then impose the weak form of the equations to improve the accuracy of the predicted solutions. 
    This method can be further generalized to any basis, time-dependent PDEs, or higher dimensional problems, by performing a change of basis, using a time discretization, or a higher dimensional convolutional kernel, respectively. The LGNet can be used to predict a wide variety of solutions to different PDEs with varying boundary conditions. %\footnote{Please see our GitHub for the code used herein. \url{https://github.com/BryceWayne/SEM} }
    %We expect our method to perform equally as well in the singular perturbation regime. 
    %The results reflected in this research will allow future in-depth studies into a wide variety of problems such as temporal, higher-dimensional or perturbative problems.

    Another interesting application for the LGNet is to study stiff partial differential equations.
    Numerical methods for convection-dominated singularly perturbed problems raise the substantial difficulties since the small diffusive parameter $\e$ produces a sharp transition inside thin layers. 
    In \cite{CHT20}, one of the authors of this paper studied an enriched spectral method to solve the singularly perturbed models. 
    It consists of adding to the Legendre basis functions, analytically-determined boundary layer elements called ``correctors", with the aim of capturing most of the complex behavior occurring near the boundary with such an element.
    Thanks to the structure of our architecture, it is easy to introduce an additional basis function with a learning parameter such that
    $u(x) \simeq \sum_{k=0}^{N-1} \alpha_k \phi_k(x) + \alpha_N \phi_N(x)$,
    where $\phi_k$, $0\leq k \leq N-1$ is a standard Legendre basis and $\phi_N$ is an enriched basis generated by the corrector function.
    Hence, we are able to extend our network to the enriched method solving singularly perturbed differential equations.
    This example manifests flexibility of our LGNet algorithm. 
    This application will be discussed in the forthcoming article.
	\newcommand{\etalchar}[1]{$^{#1}$}

%    
%    \bibliography{references}
%    \bibliographystyle{alpha}
\end{document}